# Validity of tests for time-to-event endpoints in studies with the Pocock and Simon covariate-adaptive randomization

Victoria P. Johnson[1], Michael Gekhtman[2], and Olga M. Kuznetsova[3]

Abstract

Clinical trials often have a heterogenous population with several prognostic factors. To reduce bias that can be caused by imbalances in prognostic factors between the treatment arms, randomization procedures that enforce balance in prognostic factors are commonly employed. With a small number of prognostic factors, stratified randomization is used; when the number of prognostic factors or factor levels is large, the Pocock and Simon's covariate-adaptive randomization can provide the required balance. In the presence of prognostic covariates, inference about the treatment effect with time-to-event endpoints is mostly conducted via the stratified log-rank test or the score test based on the Cox proportional hazards model. In their ground-breaking work Ye and Shao (2020) have demonstrated theoretically that when the model is misspecified, the robust score test (Wei and Lin, 1989) as well as the unstratified log-rank test are conservative in trials with stratified randomization (Ye and Shao, 2020). This fact, however, was not established for minimization other than through simulations. In this paper, we expand the results of Ye and Shao to a more general class of randomization procedures and show, in part theoretically, in part through simulations, that the Pocock and Simon covariate-adaptive allocation belongs to this class. We also advance the search for the correlation structure of the normalized within-stratum imbalances with minimization by describing the asymptotic correlation matrix for the case of equal prevalence of all strata. We expand the robust tests proposed by Ye and Shao for stratified randomization to minimization and examine their performance trough simulations.

Keywords: minimization, Pocock and Simon covariate-adaptive randomization, Type I error, log-rank test, robust survival analysis tests

[1]Corresponding author; Statistics and Data Science Innovation Hub, GlaxoSmithKline, 1250 S Collegeville Rd, Collegeville, PA 19426; Email: victoria.p.johnson@gsk.com

[2]Department of Mathematics, University of Notre Dame, Notre Dame, IN 46556

[3]Merck & Co., Inc., Kenilworth, NJ, USA

1. Introduction

In a randomized clinical trial it is desirable to have the treatment groups balanced in important prognostic factors to minimize bias in the study results. It is especially important when the prognostic factors are strong. The balance is typically established through stratified randomization that breaks the study population into strata formed by the combinations of the levels of the prognostic factors; a restricted randomization procedure (most commonly, the permuted block schedule [1]) is used within each stratum [2]. The balance within strata leads to a good marginal balance, that is balance within the levels of the factors, when the number of stratification factors is small.

However, with several strong known prognostic factors as is often the case in oncology trial, stratified randomization may become unfeasible if the population is broken into a large number of small strata [3]. In this case one of the dynamic allocation procedures that aim at marginal balance, and not the balance within individual strata, are employed [4-9]. The most frequently used and most studied dynamic allocation is the Pocock and Simon covariate-adaptive procedure [5] – a generalization of the minimization proposed by Taves [4] that adds a random element to every allocation. The Pocock and Simon procedure finds the treatment arm that would result in the smallest overall imbalance and assigns the patient to that treatment arm with high probability called bias. In practice, bias $p$ in the range of 0.8-0.95 is typically used to provide good marginal balance. In a two-arm study with equal allocation the imbalance within a level of a covariate (marginal imbalance) is typically measured by the difference in treatment totals within that level of the covariate. One common measure of the overall imbalance is the sum of the squared marginal imbalances across all levels of all factors. Another commonly used measure of the overall imbalance is the sum of the absolute marginal imbalances across all levels of all factors.

The theoretical properties of the statistical tests following the Pocock and Simon covariate-adaptive allocation were unknown for a long time and statisticians relied on simulations that established conservativeness of the tests [10-16]. However, in recent years tremendous advances into inference following covariate-adaptive minimization were made, with results on the validity of the $t$-test under linear or generalized linear model given by Shao et al. [17], Shao and Yu [18], Ma et al. [19], Bugni et al [20], and Ye [21]. However, for more complicated survival analysis tests the theoretical results were not available for either the stratified randomization or the Pocock and Simon covariate-adaptive allocation until the work of Ye and Shao [22].

In their groundbreaking work, Ye and Shao [22] developed a comprehensive theory for log-rank-type and partial likelihood score tests under stratified and covariate-adaptive randomization. They developed novel technical tools to linearize the tests and applied these tools to show that when the model is misspecified, these tests are conservative with stratified randomization. Previous results on the robustness of these tests under model misspecification (Lin and Wei [23], Kong and Slud [24], DiRienzo and Lagakos [25]) were obtained with respect to simple randomization that is rarely used in practice, therefore the results by Ye and Shao [22] signify a very important advance for clinical trials. In particular, Ye and Shao [22] theoretically demonstrated the asymptotic validity and robustness of the stratified log-rank test [26] (a commonly used analysis of a time-to-event endpoint) when all randomization strata are used in the analysis. They also proposed a variance adjustment to the partial likelihood score test that makes the test valid and robust under model misspecification.

The results of Ye and Shao [22] rely on the property of the stratified randomization that the within-stratum imbalances in treatment assignments are independent. In this paper, we expand on the methodological framework built by Ye and Shao [22] to demonstrate that the partial likelihood score test and the log-rank test are also conservative for a wider class of allocation procedures with correlated within-stratum imbalances in treatment assignments. Correlated within-stratum imbalances are common to dynamic allocation procedures, such as the Pocock and Simon covariate-adaptive procedure, for which the probability of the next treatment assignment depends on the treatment assignments of subjects already enrolled in different strata [2]. The required property of the allocation procedures in the described class is the convergence of the covariance matrix of the normalized within-stratum treatment imbalances to a non-negative-definite matrix with the maximum eigenvalue that does not exceed 1 (a part of condition 4 below). We show, in part theoretically, in part through simulations that the Pocock and Simon covariate-adaptive allocation with $p≥0.8$ belongs to the described class of the allocation procedures. As part of the research, we provide the asymptotic correlation matrix of the normalized within-stratum treatment imbalances with the Pocock and Simon covariate-adaptive allocation in the case of the equal prevalence of the strata that has not been previously described. The question of the asymptotic covariance matrix in the case of the unequal prevalence remains open.

Throughout the paper we discuss only 1:1 allocation to two treatment groups. Generalization of the allocation procedures, in particular the dynamic allocation procedures like Pocock and Simon's to unequal allocation, to unequal allocation requires special care to make sure that the allocation ratio is preserved at every allocation step [15]. Allocation procedures with unequal allocation and procedures with >2 arms are outside of the scope of this paper.

In section 2 we provide the notations, re-cap the results of Ye and Shao [22], and expand their results to a wider class of allocation procedures. In section 3 we describe the asymptotic correlation matrix for the normalized within-stratum imbalances in treatment assignments for the Pocock and Simon covariate-adaptive allocation when all strata have equal prevalence and derive its eigenvalues. We obtain the variance of the within-stratum normalized imbalances through large simulations and show that with $p≥0.8$ the product of the simulated variance and the exact maximum eigenvalue of the asymptotic correlation matrix (that is the maximum eigenvalue of the asymptotic covariance matrix) is below 1. We show through simulations that for the Pocock and Simon covariate-adaptive allocation with unequal strata prevalence, the maximum eigenvalue of the asymptotic covariance matrix is also below 1. In section 4 we apply the variance adjustment approach proposed by Ye and Shao [22] for robust tests following stratified randomization to the tests following allocation procedures that satisfy condition 4. We examine the performance of these tests following the Pocock and Simon covariate-adaptive allocation trough simulations. A discussion completes the paper.

## 2. Notations and methods

### 2.1. Notations

Consider a study with $N$ subjects randomized in 1:1 ratio to Treatment 1 and Treatment 0 in a population that has $M$ known prognostic factors (baseline covariates) $k=1, \ldots, M$; factor $k$ has $n_k >1$ levels. Denote by **Z** the vector of these baseline covariates. Thus the study population is broken into $M_s=n_1 \times \ldots \times n_M$ strata formed by the combinations of the factor levels; these strata correspond to the categories of **Z**.

Denote the stratum $z$ as $z=(i_1, \ldots, i_M)$, where $i_1, \ldots, i_M$ are the levels of the factors that formed the stratum $z$. For example, if age (<65, ≥65) and family history of disease (Yes/No) are two binary prognostic factors and smoking status is the third factor that has three level (Never smoke/ former smoker/ current smoker), $M=3$, $n_1=n_2=2$, $n_3=3$, $M_s=12$ and the stratum "≥65, No family history, never smoke" is denoted by $z=(2,2,1)$. We will also use the notation where the strata are numbered from 1 to $M_s$.

Denote by $Z_i=(Z_{i1}, \ldots, Z_{iM})$, $i=1, \ldots, N$, the vector of baseline covariates of the $i$-th subject; this vector determines the stratum the $i$-th subject belong.

Denote by $N_z$ the size of the stratum $z$; $\sum_{z=1}^{M_s} N_z = N$. Let $w_z = N_z/N$ be the fraction of subjects in the stratum $z$, $\sum_{z=1}^{M_s} w_z = 1$.

Denote by $m_z$ the size of the stratum $z$ after $m$ allocations and by $D_m(z)$ the imbalance in treatment assignments within stratum $z$, that is the difference in the number of subjects allocated to Treatment 1 and Treatment 0, after $m$ allocations. Let us denote, and by $d(z)$ the imbalance within stratum $z$ normalized by its size (normalized imbalance) after $m$ allocations: $d_m(z) = \frac{D_m(z)}{\sqrt{m_z}}$. Then

$$d_N(z) = \frac{D_N(z)}{\sqrt{N_z}} = \frac{D_N(z)}{\sqrt{N}\sqrt{w_z}} \tag{1}$$

Denote by $D_m(j; h)$ the marginal imbalance in treatment assignments, that is the difference in the number of subjects allocated to Treatment 1 and Treatment 0 among the subjects with $i_j=h$, after $m$ allocations. Then

$$D_m(j; h) = \sum_{z: i_j = h} D_m(z) \tag{2}$$

For example, the marginal imbalance across smokers is equal to the sum of the within-stratum imbalances across the 4 strata that include the smokers. When the marginal imbalance $D_N(j; h)$ is low, the two treatment arms have very similar percentages of subjects with $i_j=h$. Typically, it is desired to get similar distributions of the levels of a strong prognostic factor in each group, that is low marginal imbalances, while the within-strata balance is of lesser importance. However, the only non-dynamic way of enforcing marginal balance is to use stratified randomization that breaks the study population into strata formed by the combinations of factor levels; a restricted randomization, for example, permuted block randomization, is then used to randomize subjects within each stratum. Balance in treatment assignments within strata then enforces the marginal balance when the number of strata is not too large.

A dynamic Pocock and Simon covariate-adaptive allocation aims to minimize the marginal imbalances but not the imbalances within individual strata $D_N(z)$. We will use the classical Pocock and Simon covariate-adaptive procedure in this paper where the imbalance within a level of a covariate after allocation of *m* subjects is measured by the difference in treatment totals across all subjects with this level of the covariate and the overall imbalance is calculated as the sum of the squared imbalances across all levels of all covariates:

$$Imb(m) = \sum_{j=1}^{M} \sum_{h=1}^{n_j} [D_m(j;h)]^2 \qquad (3)$$

When a new subject arrives for randomization, the treatment arm that would result in the lowest overall imbalance (3) after the allocation is assigned to the patient with probability $p>1/2$ (we will use bias $p=0.9$ in most simulations). If both arms would result in the same overall imbalance, Treatment 1 is assigned with probability 0.5.

If it is important to have balance in treatment assignments within certain individual strata, these strata can be added to the minimization scheme as separate factors and thus balanced [5]. Hu and Hu [27] and Hu and Zhang [28] consider a generalization of the Pocock and Simon minimization where the treatment arm imbalance in the number of subjects, all marginal imbalances and all within-stratum imbalances are given non-negative weights. For the classical Pocock and Simon covariate adaptive procedure we consider in this paper the weights assigned to imbalances within individual strata and the imbalance in the treatment arm totals are 0.

Denote by $Cov(\boldsymbol{d}_N(z)) = [cov(d_N(i_1,…,i_M), d_N(j_1,…,j_M))]$, $i_1=1, …, n_1$; ….; $i_M=1, …, n_M$; $j_1=1, …, n_1$; ….; $j_M=1, …, n_M$, the covariance matrix of the vector of normalized imbalances $\boldsymbol{d}_N(z)$, $z=1, …, M_s$, and by $Cor(\boldsymbol{d}_N(z)) = [cor(d_N(i_1,…,i_M), d_N(j_1,…,j_M))]$, $i_1=1, …, n_1$; ….; $i_M=1, …, n_M$; $j_1=1, …, n_1$; ….; $j_M=1, …, n_M$, the correlation matrix of the vector of normalized imbalances $\boldsymbol{d}_N(z)$, $z=1, …, M_s$.

The normalized within-stratum imbalances and their correlation depend on the procedure used to randomize the study subjects. For the stratified permuted block randomization, a commonly used in clinical trial randomization procedure, an imbalance within a stratum does not exceed half of the block size, thus the normalized imbalance (1) converges to 0 as $N_z$ increases and $Cov(\boldsymbol{d}_N(z)) \to O$, where $O$ is the $M_s \times M_s$ matrix with all terms equal to 0. Similar asymptotic behavior holds when one of the allocation procedures with maximum tolerated imbalance that require $|D_N(z)|<b$ for prespecified $b>0$ [29], such as maximal procedure by Berger, Ivanova Knoll [30], big stick design by Soares and Wu [31], biased coin design with imbalance tolerance by Chen [32] is used to randomize subjects within each stratum. For stratified randomization based on the biased coin randomization by Efron [33], where a subject is assigned with probability $p>1/2$ to the treatment arm that has fewer subjects within a stratum and in 1:1 ratio when both arms have the same number of subjects, $D_N(z)$ the is bounded in probability [33] and thus $Cov(\boldsymbol{d}_N(z)) \to O$. For a stratified randomization based on the urn design by Wei [34,35], $Cov(\boldsymbol{d}_N(z)) \to \nu I$ where $\nu = 1/3$ based on [34, 35] and $I$ is the $M_s \times M_s$ identity matrix.

For complete randomization, $Cov(\boldsymbol{d}_N(z)) \to I$. It should be noted that in practice complete randomization is not used with stratified randomization as it would be the same as applying complete randomization without stratification and will not promote the balance in covariates. However, Ye and Shao [22] include stratified complete randomization in their paper for reference.

While for stratified randomization the within-stratum imbalances in two different strata are independent, this is not the case for the Pocock and Simon covariate-adaptive allocation. The covariance matrix of the within-strata imbalances was not described for the Pocock and Simon allocation, but Hu and Zhang [28] showed that when individual strata are assigned 0 weight in the calculation of the overall imbalance (as in the version considered in this paper), for every stratum $z$

$$d_N(z) = \frac{D_N(z)}{\sqrt{N_z}} \xrightarrow{D} N(0, \sigma_z) \tag{4}$$

where $\xrightarrow{D}$ denotes convergence in distribution as $N \to \infty$.

This implies that the diagonal elements of the covariance matrix of the normalized imbalances $Cov(\boldsymbol{d}_N(z))$ converge to $\sigma_z$ as $N \to \infty$.

Hu and Zhang [28] also showed that the marginal imbalances $D_m(j; h)$ are bounded in probability for all $j=1,...,M$, $h=1,..., n_j$. Therefore, from (2)

$$\sum_{z: i_j=h} \frac{D_N(z)}{\sqrt{N}} \xrightarrow{P} 0. \tag{5}$$

Hu and Zhang demonstrated these properties for the Pocock and Simon covariate-adaptive allocation with the overall imbalance measured as the sum (3) of the squared imbalances across all levels of all covariates. However, simulations suggest that the same property holds for other measures of the overall imbalance, for example, the sum of the absolute imbalances across all levels of all covariates.

Ye and Shao [22] introduce condition 1 common to many allocation procedures.

<u>Condition 1</u>. $\boldsymbol{Z}$ is discrete with finitely many categories, $E(I_i|\boldsymbol{Z}_1, ..., \boldsymbol{Z}_N) = \frac{1}{2}$ for $i=1, ..., N$ and

$$\frac{D_N(z)}{N_z} = o_p(1) \text{ for every category } z \text{ of } \boldsymbol{Z}. \tag{6}$$

The condition $E(I_i|\boldsymbol{Z}_1, ..., \boldsymbol{Z}_N) = \frac{1}{2}$ holds for all equal allocation procedures with randomization algorithm symmetric with respect to the two treatment arms. Although an equal allocation procedure with an asymmetric randomization algorithm can be easily invented, such procedures are rarely used in clinical trials. As shown in Baldi et al. [36], (6) holds for most randomization procedures with equal allocation, in particular, it obviously holds for all 1:1 allocation procedures mentioned in this paper (for the Pocock and Simon procedure (6) follows from (4)).

They also introduced conditions 2 and 3 that describe the two classes of randomization procedures for which they prove their results.

Condition 2. $\frac{D_N(z)}{\sqrt{N_z}} = o_p(1)$ for every $z$.

Condition 3. $\frac{D_N(z)}{\sqrt{N_z}} \xrightarrow{D} N(0, v)$, where $v \leq 1$, and $D_N(z_1)$ and $D_N(z_2)$ are independent for all $z_1 \neq z_2$.

While all stratified randomization procedures mentioned above satisfy either condition 2 or 3, the Pocock and Simon covariate-adaptive procedure does not satisfy either of these conditions.

We will show that the results of Ye and Shao [22] can be expanded to a class of randomization procedures that satisfy condition 4 - a more general condition than condition 3.

Condition 4. Property (4) holds and $Cov(\boldsymbol{d}_N(z)) \xrightarrow{N \to \infty} Cov$, where $Cov$ is a non-negative definite $M_s \times M_s$ matrix with the maximum eigenvalue $\Lambda_{max} \leq 1$.

For the allocation procedures that satisfy condition 3, $Cov = vI$, $v \leq 1$; therefore the maximum eigenvalue of $Cov$ is $\Lambda_{max} = v$. Thus, condition 3 implies condition 4.

For the Pocock and Simon covariate-adaptive randomization the covariance matrix of the normalized within-strata imbalances is not diagonal as the imbalances in different strata are correlated, but (4) is satisfied [26]. The asymptotic covariance matrix $Cov$ has not been described in the literature; we will show in Section 3– in part theoretically, in part through simulations – that for bias $p \geq 0.8$ its maximum eigenvalue is $\leq 1$ and thus, the Pocock and Simon procedure satisfies condition 4.

2.2. Methods

Ye and Shao [22] developed tools for studying the score tests in trials with allocation procedures that satisfy condition 1 on the vector of baseline covariates $\boldsymbol{Z}$ for which we want to balance the randomization.

We will follow the survival analysis notations and concepts in Ye and Shao [22] paper that we will repeat below.

Let $\boldsymbol{V}_i$ be the vector of all measured and unmeasured covariates, and let $I_i$ be the treatment indicator (0 or 1) of the $i$-th subject, $i=1, \ldots, N$. Let $X_j^*$ and $C_j$ be the potential failure time and the censoring time for a patient assigned treatment $j$, $j=0,1$. Let us assume that $(X_{i0}^*, C_{i0}, X_{i1}^*, C_{i1}, \boldsymbol{V}_i)$, $i=1, \ldots, N$, are independent and identically distributed samples from $(X_0^*, C_0, X_1^*, C_1, \boldsymbol{V})$. Then the observed response with right censoring for subject $i$ is $(X_i, \delta_i)$, where $X_i = I_i X_{i1} + (1 - I_i) X_{i0}$, $X_{ij} = min(X_{ij}^*, C_{ij})$, $\delta_i = I_i \delta_{i1} + (1 - I_i) \delta_{i0}$, $\delta_{ij} = I(X_{ij} = X_{ij}^*)$, and $I(\cdot)$ is the indicator function.

Ye and Shao [22] make 3 assumptions on the censoring distribution which hold under many realistic situations in clinical trials (see more on that in [22]).

*Assumption 1.* $(X_{ij}^*, C_{ij}, \boldsymbol{V}i)s$ and $I_i s$ are conditionally independent given $\boldsymbol{Z}_i s$.

*Assumption 2.* $X_j^*$ and $C_j$ are independent conditioned on covariate $\boldsymbol{V}$, $j=0,1$.

*Assumption 3.* The ratio $P(C_1 \geq t | \boldsymbol{V}) / P(C_0 \geq t | \boldsymbol{V})$ is a function of $t$ only.

Consider the working Cox proportional hazard model [37, 38] where the true hazard of the *i*-th subject $\lambda(t, \boldsymbol{V}_i, I_i)$ is modelled through a vector of covariates $\boldsymbol{W}_i$, vector of unknown parameters $\beta$, an unknown treatment effect parameter $\theta$ and a baseline hazard function $\lambda_0(t)$:

$$\lambda(t, \boldsymbol{V}_i, I_i) = \lambda_0(t) exp(\theta I_i + \beta' \boldsymbol{W}_i). \tag{7}$$

The partial likelihood function based on working model above is

$$L(\theta, \beta) = \prod_{i=1}^{N} \left( \frac{exp(\theta I_i + \beta' \boldsymbol{W}_i)}{\sum_{k=1}^{N} Y_k(X_i) exp(\theta I_i + \beta' \boldsymbol{W}_i)} \right)^{\delta_i},$$

where

$$Y_i(t) = I_i Y_{i1}(t) + (1 - I_i) Y_{i0}(t), \ Y_{ij}(t) = I(X_{ij} \geq t). \tag{8}$$

Let $\hat{\beta}_0$ be the maximum partial likelihood estimator of $\beta$ under the constraint $\theta = 0$. From theorem 2.1 of Struthers and Kalbfleisch [39], $\hat{\beta}_0$ converges in probability to a unique vector $\beta_*$ under some regularity conditions (see the conditions in Theorem 1) regardless of whether the working model above is correct or not, where $\beta_* = \beta$ if the working model is true.

Ye and Shao [22] examined the asymptotic behavior of the score test $T_S$ by Lin and Wei [23]:

$$T_S = N^{-\frac{1}{2}} U_\theta(0, \hat{\beta}_0) / \hat{B}(0, \hat{\beta}_0)^{1/2} \tag{9}$$

where

$$U_\theta(0, \beta) = \frac{\partial log(L(\theta, \beta))}{\partial \theta} \bigg|_{\theta=0} = \sum_{i=1}^{N} \int_0^\tau \left\{ I_i - \frac{S_N^{(1)}(\beta, t)}{S_N^{(0)}(\beta, t)} \right\} dN_i(t)$$

$$S_N^{(r)}(\beta, t) = N^{-1} \sum_{i=1}^{N} Y_i(t) exp(\beta' \boldsymbol{W}_i) I_i^r, r = 0,1,$$

$$N_i(t) = I_i N_{i1}(t) + (1 - I_i) N_{i0}(t), \ N_{ij}(t) = \delta_{ij} I(X_{ij} \leq t),$$

$\tau$ is a point satisfying $P(X_{ij} \geq \tau) > 0$ for *j*=0,1, and

$$\hat{B}(0, \hat{\beta}_0) = N^{-1} \sum_{i=1}^{N} \left[ \delta_i \left\{ I_i - \frac{S_N^{(1)}(\beta, X_i)}{S_N^{(0)}(\beta, X_i)} \right\} - \sum_{j=1}^{N} \frac{\delta_i Y_i(X_j) exp(\beta' \boldsymbol{W}_i)}{N S_N^{(0)}(\beta, X_j)} \left\{ I_i - \frac{S_N^{(1)}(\beta, X_j)}{S_N^{(0)}(\beta, X_j)} \right\} \right]^2$$

is a robust variance estimator under simple randomization in the sense that it is consistent regardless of whether the working model (7) is correct or not.

They also examined the unstratified log-rank test, a special case when $\boldsymbol{W}_i \equiv 0$, that is, no covariates are used in the working model. In this case $N^{-\frac{1}{2}} U_\theta(0, \hat{\beta}_0)$ in (9) equals the numerator of the log-rank test statistic

$$T_L = N^{-\frac{1}{2}} \sum_{i=1}^{N} \int_0^\tau \left\{ I_i - \frac{\bar{Y}_1(t)}{\bar{Y}(t)} \right\} dN_i(t) / \hat{\sigma}, \tag{10}$$

where $\bar{Y}_1(t) = \sum_{i=1}^{N} I_i Y_i(t), \bar{Y}_0(t) = \sum_{i=1}^{N} (1 - I_i) Y_i(t), \bar{Y}(t) = \bar{Y}_1(t) + \bar{Y}_0(t)$, and

$$\hat{\sigma}^2 = \frac{1}{N}\sum_{i=1}^{N}\int_0^\tau \frac{\bar{Y}_1(t)\bar{Y}_0(t)}{\bar{Y}^2(t)}dN_i(t) \ . \tag{11}$$

The following notations from Ye and Shao [22] are required to describes their results; the expectation $E$ and variance $var$ are taken under $H_0$.

Let

$$\mu(t) = E\{I_i|Y_i(t) = 1\}$$

$$p(t) = E\{Y_i(t)\lambda(t,V_i)\}/E\{Y_i(t)exp(\beta_*'\boldsymbol{W}_i)\}$$

$$O_{i1} = \int_0^\tau \{1 - \mu(t)\}\{dN_{i1}(t) - p(t)Y_{i1}(t)exp(\beta_*'\boldsymbol{W}_i)dt\} \tag{12}$$

$$O_{i0} = \int_0^\tau \{1 - \mu(t)\}\{dN_{i0}(t) - p(t)Y_{i0}(t)exp(\beta_*'\boldsymbol{W}_i)dt\} \tag{13}$$

Per Ye and Shao [22] for 1:1 allocation procedures that satisfy condition 1 and under assumptions 1-3 on censoring distribution, the numerator of the score test statistics can be expressed as

$$N^{-1/2}U_\theta(0,\hat{\beta}_0) = U_1 + U_2 + o_p(1),$$

where

$$U_1 = \frac{1}{\sqrt{N}}\sum_{i=1}^{N}\{I_i(O_{i1} - E_i) - (1 - I_i)(O_{i0} - E_i)\},$$

$$U_2 = \frac{1}{\sqrt{N}}\sum_{i=1}^{N}(2I_i - 1)E_i, \tag{14}$$

$$E_i = E(O_{ij}|\boldsymbol{Z}_i). \tag{15}$$

As Ye and Shao [22] demonstrated, under $H_0$, condition 1, and under assumptions 1-3, $E_i$ in (15) does not depend on $j$,

$$E(E_i) = 0, \tag{16}$$

$E(U_1) = E(U_2) = 0$, and

$$U_1 \xrightarrow{D} N(0,\psi) \ ,$$

where

$$\psi = E\{var(O_{i1}|\boldsymbol{Z}_i) + var(O_{i0}|\boldsymbol{Z}_i)\}/2 \ . \tag{17}$$

Let us denote

$$\phi = var(E_i) \ .$$

From (15) and (16)

$$\phi = \sum_z w_z E^2(O_{ij}|\boldsymbol{Z}_i = z). \tag{18}$$

Consider vector $\boldsymbol{G} = (G_z, z = 1,\ldots,M_s)'$, where $G_z = \sqrt{w_z}E(O_{ij}|\boldsymbol{Z}_i = z)$. $\tag{19}$

Note that

$$\boldsymbol{G}'\boldsymbol{G} = var(E_i) = \phi.$$

From (14) and (1) $U_2$ can be written as

$$U_2 = \frac{1}{\sqrt{N}}\sum_z D_n(z)E(O_{ij}|\boldsymbol{Z}_i = z) = \sum_z d_N(z)\sqrt{w_z}E(O_{ij}|\boldsymbol{Z}_i = z) = \sum_z d_N(z)G_z \quad (20)$$

Thus, from (20), the variance of $U_2$ can be written as

$$Var(U_2) = \boldsymbol{G}'Cov(\boldsymbol{d}_N(z))\boldsymbol{G}$$

Note that if $\phi = 0$, then $E(O_{ij}|\boldsymbol{Z}_i = z)$ are the same for all $z$ and $Var(U_2) = 0$.

Under condition 4,

$$Var(U_2) \xrightarrow{N\to\infty} \boldsymbol{G}'Cov\boldsymbol{G}$$

also, since $\Lambda_{max} \leq 1$,

$$\boldsymbol{G}'Cov\boldsymbol{G} \leq \boldsymbol{G}'\boldsymbol{G} = \phi$$

Ye and Shao also demonstrated that under $H_0$, condition 1, and assumptions 1-3, $U_1$ and $U_2$ are uncorrelated, and that the denominator of the Lin and Wei [23] score test statistics converges in probability to

$$\hat{B}(0,\hat{\boldsymbol{\beta}}_0) \xrightarrow{P} \psi + \phi.$$

While Ye and Shao built the foundation for deriving the asymptotic properties of the survival tests under covariate-adaptive randomization, the theoretical part of their break-through paper was limited to the examination of the stratified randomization procedures for which $Cov = vI, 0 \leq v \leq 1$. However, their methodology can be expanded to allocation procedures that satisfy a more general condition 4. In particular, Theorem 1 in Ye and Shao can be expanded by replacing the asymptotic variance of $N^{-\frac{1}{2}}U_\theta(0,\hat{\beta}_0)$ with $\psi + \boldsymbol{G}'Cov\boldsymbol{G}$ that is appropriate under condition 4.

*Theorem 1.* Consider the partial likelihood function under working model (7). Assume condition 1 and censoring assumptions 1-3 and that $-N^{-1}\partial^2 log\{L(\theta,\beta)\}/\partial\beta\partial\beta'|_{\beta=\beta}$ converges in probability to a positive-definite matrix. Then under $H_0$ and randomization that satisfies condition 4,

$$N^{-\frac{1}{2}}U_\theta(0,\hat{\beta}_0) \xrightarrow{D} N(0,\psi + \boldsymbol{G}'Cov\boldsymbol{G})$$

where $\psi$ and $\boldsymbol{G}$ are defined in (17) and (19).

Similarly, the corollary 1 by Ye and Shao [22] that establishes the asymptotic behavior of the score test statistics $T_S$ for stratified randomization procedures that satisfy conditions 1 and condition 2 or 3, can be re-stated for a class of randomization procedures that satisfy condition 4.

*Corollary 1.* Under the same assumptions as in Theorem 1, $T_S$ in expression (9) has the following property under $H_0$, regardless of whether the working model (7) is correct or not:

$$\lim_{N\to\infty} P(|T_S| > Z_{\frac{\alpha}{2}}) = 2\Phi\left\{-Z_{\frac{\alpha}{2}}\left(\frac{\psi+\phi}{\psi+G'[Cov]G}\right)^{1/2}\right\}, \quad (21)$$

where $\Phi$ is a standard normal cumulative distribution function.

Per Ye and Shao [22], if $\phi = 0$ (as is the case, in particular, when the working model (7) is correct), $U_2 = o_p(1)$ and both the numerator and denominator in (21) are equal to $\psi$ which means that $T_S$ is valid.

If $\phi > 0$, let us denote

$$\gamma = G'CovG/\phi = G'CovG/G'G \quad (22)$$

In this case from (21) the variance of the asymptotic distribution of $T_s$ can be written as

$$\frac{\psi + G'CovG}{\psi + \phi} = \frac{\psi + \gamma\phi}{\psi + \phi} \quad (23)$$

*Corollary 1.1.* Under condition 4, $T_S$ is either valid or conservative. If $\phi = 0$ or $\gamma = 1$ for vector $G$, then the test is valid; otherwise the test is conservative.

Corollary 1.1 follows directly from Corollary 1, (23), and the fact that $\Lambda_{max} \leq 1$ in condition 4 implies that for any vector $G$ such that $\phi > 0$, $\gamma \leq 1$.

Similarly to corollary 1, corollary 2 by Ye and Shao [22] that establishes the asymptotic behavior of the unstratified log-rank test statistics $T_L$ can be re-stated for a class of randomization procedures that satisfy condition 4. Ye and Shao demonstrated that under $H_0$ the variance estimator in (11)

$$\hat{\sigma}^2 \xrightarrow{P} \tilde{\psi} + \tilde{\phi}$$

where $\tilde{\psi} = E\{var(\tilde{O}_{i1}|Z_i) + var(\tilde{O}_{i0}|Z_i)\}/2$

and $\tilde{\phi} = var\{E(\tilde{O}_{ij}|Z_i)\}$,

which are special cases of $\psi$ and $\phi$ defined in (17) and (18) by setting $W_i \equiv 0$ in $O_{ij}$ defined in (12) and (13).

Defining vector $\tilde{G} = (\sqrt{w_z}E(\tilde{O}_{ij}|Z_i = z), z = 1, \ldots, M_s)'$,

and noting that $\tilde{G}'\tilde{G} = var\{E(\tilde{O}_{ij}|Z_i)\} = \tilde{\phi}$, we can re-state corollary 2 of Ye and Shao [22] as following:

*Corollary 2.* Under condition 1 and assumptions 1-3, $T_L$ in expression (10) has the following property under $H_0$ and randomization that satisfies condition 4

$$\lim_{N\to\infty} P(|T_L| > Z_{\frac{\alpha}{2}}) = 2\Phi\left\{-Z_{\frac{\alpha}{2}}\left(\frac{\tilde{\psi}+\tilde{\phi}}{\tilde{\psi}+\tilde{G}'Cov\tilde{G}}\right)^{1/2}\right\} \quad (24)$$

Similar to the score test, per [22], if $\tilde{\phi} = 0$, $U_2=o_p(1)$ and both the numerator and denominator in (24) are equal to $\tilde{\psi}$ which means that $T_L$ is valid. Per [22], $\tilde{\phi} = 0$ when outcome is independent of $Z$.

For $\tilde{\phi} > 0$ let us define

$$\tilde{\gamma} = \frac{\tilde{G}' Cov \tilde{G}}{\tilde{\phi}}. \tag{25}$$

*Corollary 2.1.* Under condition 4, $T_L$ is either valid or conservative. If $\tilde{\phi} = 0$ or $\tilde{\gamma} = 1$ for vector $\tilde{G}$, then the log-rank test is valid; otherwise the log-rank test is conservative.

Corollary 2.1 follows directly from Corollary 2, (23), and the fact that $\Lambda_{max} \leq 1$ in condition 4 implies that for any vector $\tilde{G}$ such that $\tilde{\phi} > 0$, $\tilde{\gamma} \leq 1$.

In Section 3 we will show (in part theoretically, in part through simulations) that the Pocock and Simon covariate-adaptive randomization with bias $p \geq 0.8$ satisfies condition 4 and thus, the asymptotic behavior of the score test and the log-rank test following the Pocock and Simon randomization is described by corollary 1.1 and 2.1, respectively.

Note that while for a stratified randomization with $Cov = vI$ the product $G'CovG$ and thus $\gamma$ does not depend on vector $G$, for other non-negative matrices $Cov$, $\gamma$ in (22) depends on $G$. Thus the degree of conservativeness of $T_S$ when $Cov \neq vI$ depends on vector $G$. Similarly, $\tilde{\gamma}$ in (25) depends on vector $\tilde{G}$, and thus the degree of conservativeness of $T_L$ for matrices $Cov \neq vI$ depends on vector $\tilde{G}$.

For example, as will be demonstrated in the Appendix 2, for the Pocock and Simon covariate-adaptive allocation with bias $p=0.9$ that balances on two factors with 2 levels each where all 4 strata have the same prevalence, the asymptotic covariance matrix $Cov$ has a single positive eigenvalue (0.94) and 3 other eigenvalues are equal to 0. The eigenvector that corresponds to the maximum eigenvalue is $H$=(1, -1, -1, 1). If $G$ is close to being collinear with $H$, then $\gamma = G'CovG/G'G$ is closer to 1 (the test is less conservative); if $G$ is close to being orthogonal to $H$, then $\gamma = G'CovG/G'G$ is close to 0 (the test is more conservative).

3. Covariance and correlation matrices for the normalized within-stratum imbalances for the Pocock and Simon covariate-adaptive allocation

Since for a single factor ($M=1$) the Pocock and Simon covariate-adaptive procedure degenerates into a stratified randomization with biased coin randomization performed within each stratum, we will consider only the case of $M \geq 2$.

Consider first the case when all strata $z$ have the same prevalence, that is for large $N$ for all $z$

$$w_z = \frac{1}{n_1 \times ... \times n_M} \text{ and}$$

$$N_z = \frac{N}{n_1 \times ... \times n_M}. \tag{26}$$

This implies that all levels of any factor $k$ have the same prevalence ($1/n_k$) and thus, the Pocock and Simon covariate-adaptive procedure is symmetric with respect to all levels of factor $k$. In this case the covariance matrix $Cov(\mathbf{d}(z))$ is invariant with respect to any permutation of the levels of any factor, that is, for any set of permutations $f_1, \ldots, f_M$ of the levels of the $1^{st}, \ldots, M^{th}$ factor (including the permutation that leaves all levels in place),

$$cov(d_N(i_1, \ldots, i_M), d_N(j_1, \ldots, j_M)) = cov(d_N(f_1(i_1), \ldots, f_M(i_M)), d_N(f_1(j_1), \ldots, f_M(j_M))) \quad (27)$$

For example, when randomization is balanced on two factors with two levels each,

$$cov(d_N(1,2), d_N(1,1)) = cov(d_N(2,2), d_N(2,1)) = cov(d_N(1,1), d_N(1,2)) = cov(d_N(2,1), d_N(2,2))$$

It follows from (27) that the asymptotic correlation matrix $Cov$ is also invariant with respect to any permutation of the levels of any factor and also that the variance $\sigma_z$ defined in (4) is the same all $z$, and thus

$$Cov = \sigma_z^2 Cor \quad (28)$$

Therefore, the matrix $Cor$ is also invariant with respect to any permutation of the levels of any factor.

The following notations will be used to describe the asymptotic correlation matrix $Cor$ and its eigenvalues.

For two strata $(i_1, \ldots, i_M)$ and $(j_1, \ldots, j_M)$ let $\varepsilon_l$, $l=1, \ldots, M$, be the indicator of the two strata having the same level of the factor $l$:

$\varepsilon_l=1$ if $i_l=j_l$

$\varepsilon_l=0$ if $i_l \neq j_l$.

For example, consider minimization that balances on age (<65, ≥65), family history of disease (Yes/No), and smoking status (Never smoke/former smoker/current smoker). For the two strata $z_1=(2,2,1)$ ("≥65, No family history, never smoke") and $z_1=(1,2,1)$ ("<65, No family history, never smoke"), $\varepsilon_l=0$, $\varepsilon_2=1$, $\varepsilon_3=1$.

Since the matrix $Cor$ is invariant with respect to any permutation of the levels of any factor, the correlation coefficient that corresponds to two strata $(i_1, \ldots, i_M)$ and $(j_1, \ldots, j_M)$ depends only on the vector $\varepsilon_1, \ldots, \varepsilon_M$ and thus can be written as

$$cor((i_1, \ldots, i_M), (j_1, \ldots, j_M)) = c_{\varepsilon_1, \ldots, \varepsilon_M}.$$

The following more convenient notation for $c_{\varepsilon_1, \ldots, \varepsilon_M}$ will also be used: for a subset $I$ of $\{1, \ldots, M\}$, denote

$c_I = c_{\varepsilon_1, \ldots, \varepsilon_M}$, where $\varepsilon_l=1$ if $l \in I$ and $\varepsilon_l=0$ otherwise. In our 3-factor example ($M=3$), $c_\emptyset = c_{000}$; $c_{\{1\}} = c_{100}$; $c_{\{2\}} = c_{010}$; $c_{\{3\}} = c_{001}$; $c_{\{1,2\}} = c_{110}$; $c_{\{1,3\}} = c_{101}$; $c_{\{2,3\}} = c_{011}$; $c_{\{1,2,3\}} = c_{111}$.

Thus, the terms of the asymptotic correlation matrix *Cor* will be written as $cor((i_1, \ldots, i_M), (j_1, \ldots, j_M)) = c_I, I \subset \{1, \ldots, M\}$, where $I$ denotes the set of the factors that have the same levels for the two strata $(i_1, \ldots, i_M)$ and $(j_1, \ldots, j_M)$. In these notations $c_{\{1,\ldots,M\}} = 1$.

Denote by $I^C$ the complement of $I$ in $\{1, \ldots, M\}$: $I^C = \{1, \ldots, M\} \backslash I$; denote by $\#I$ the number of elements of $I$.

Note that in the case of the equal strata prevalence, from (26) and (5), for any level $h$ of factor $i_j$ the sum of the normalized within-stratum imbalances over all strata $z$ such that $i_j = h$ converges to 0:

$$\sum_{z:i_j=h} d_N(z) = \sum_{z:i_j=h} \frac{D_N(z)}{\sqrt{N_z}} \xrightarrow{P} 0 \tag{29}$$

Property (29) and the fact that the matrix *Cor* is invariant with respect to any permutation of the levels of any factor allows one to write a system of linear equations on its terms $c_I$.

For factor $j$, consider the set of strata $\{z: i_j = 1\}$. There are two types of relationship a stratum can have with this set: it can either belong to this set (that is have $i_j = 1$) or not belong to this set (have $i_j \neq 1$). Consider two such strata: stratum $z_1$: $i_l=1$, $l=1,\ldots, k$ that belongs to this set and stratum $z_2$: $i_j=2$, $i_l=1$, $l \neq j$, $l=1,\ldots, k$ that does not belong to the set $\{z: i_j = 1\}$.

It follows from (29) and (4) that for any stratum $z_0$, the covariance of the sum of the normalized within-stratum imbalances over all strata $z$ such that $i_j =1$ and $d_N(z_0)$ converges to 0; in particular, this takes place for $z_0 = z_1$ and $z_0 = z_2$:

$$cov\left(\sum_{z:i_j=1} d_N(z), d_N(z_1)\right) \xrightarrow{N \to \infty} 0$$

$$cov\left(\sum_{z:i_j=1} d_N(z), d_N(z_2)\right) \xrightarrow{N \to \infty} 0$$

Therefore,

$$\sum_{z:i_j=1} cov(z, z_1) = 0 \tag{30}$$

$$\sum_{z:i_j=1} cov(z, z_2) = 0 \tag{31}$$

Considering (28), and accepting a convention $\prod_{i \in \emptyset} R_i = 1$, equations (30)-(31) can be written as

$$\sum_{I \subset \{1,\ldots,M\}, j \in I} c_I \prod_{i \in I^C} (n_i - 1) = 0, j=1, \ldots, k \tag{32}$$

$$\sum_{I \subset \{1,\ldots,M\}, j \notin I} c_I \prod_{i \in I^C \backslash \{j\}} (n_i - 1) = 0, j=1, \ldots, k \tag{33}$$

Indeed, when $j \in I$, there are $\prod_{i \in I^C}(n_i - 1)$ strata among $z: i_j = 1$ that have common levels with $z_1$ for all factors in $I$ and only these factors, that is have level 1 for all factors in $I$ and level other than 1 for all factors that do not belong to $I$. When $j \notin I$, there are $\prod_{i \in I^C \backslash \{j\}}(n_i - 1)$ strata among $z: i_j = 1$ that have common levels with $z_2$ for all factors in $I$ and only these factors, that is have level 1 for all factors in $I$, and level other than 1 for all factors in $\{1, \ldots, M\} \backslash \{j\} \backslash I$.

The system of $2M$ equations (32)-(33) has $2^M -1$ independent variables $c_I, I \subset \{1, ..., M\}$, $I \neq \{1, ..., M\}$. The system (32)-(33) is not of full rank. When $M=2$, there are 3 independent equations and 3 independent variables $c_\emptyset, c_{\{1\}}$, and $c_{\{2\}}$, thus the system (32)-(33) (that is easy to solve in this case) has a unique solution:

$$c_\emptyset = \frac{1}{(n_1-1)(n_2-1)}, \; c_{\{1\}} = \frac{1}{(n_2-1)}, \; c_{\{2\}} = \frac{1}{(n_1-1)}. \tag{34}$$

When $M>2$, the solution is not unique. Nevertheless, we observe the following:

Statement. When all strata have equal prevalence,

$c_I = 1$, when $I=\{1, ...., M\}$;

Otherwise,

$$c_I = \frac{M-1-\sum_{i \in I} n_i}{\prod_{i=1}^{M} n_i - \sum_{i=1}^{M} n_i + M-1} \tag{35}$$

This statement was derived empirically by examining the correlation coefficients obtained from the simulations for $M>2$; for $M=2$, (35) becomes (34). We will consider it a true statement even though we do not have a theoretical proof as for any specified set of $M$ and $n_1, ..., n_M$ one can perform the simulations and confirm (35). It is unclear what additional requirements beyond (32)-(33) should be met for a correlation matrix to arrive at this solution; this question needs to be further explored. The equations (32)-(33) do not depend on the bias $p$ and the simulations conducted for different values of $p$ resulted in the same correlations matrix (35) (results are provided for $p=0.9$ only). We conclude from this that the asymptotic correlation matrix of normalized imbalances $Cor$ does not depend on the bias $p$.

Table A1 in Appendix 1 provides the comparison of the correlation coefficients derived through simulations with the correlation coefficients in (35) for 4-factor examples of equal strata prevalence minimization with $p=0.9$. The simulations included generation of 10,000 random sets of independent vectors of 4 independent baseline covariates with equal distribution of the levels of each factor that were allocated using the Pocock and Simon covariate-adaptive randomization. The size of a simulated dataset was $N=500 \times M_s$, so that on average each stratum would have 500 subjects. The covariance matrix of the normalized within-stratum imbalances at the end of randomization was estimated across all simulations. The variance of the normalized within-stratum imbalances $\sigma_z^2$ was estimated by pooling imbalances across all strata. For every $I \subset \{1, ..., M\}$ the estimate of $c_I$ was derived as an average of the estimated correlations across pairs of strata that have common levels of all factors in $I$ and no other factors.

The simulations presented in Table A1 were performed using the sum of the squared imbalances across all levels of all covariates as the measure of the overall imbalance (as in Hu and Zhang [28]). However, we obtained the same results when using the sum of the absolute imbalances across all levels of all covariates as the measure of the overall imbalance (results not provided).

Some observations based on (35) can be offered. Note that the denominator of (35) does not depend on $I$ and is always positive; it can be interpreted as the number of independent equations in the system (32)-(33). When $I = \emptyset$, that is the two strata have no factors in common, the numerator is $M$-1 and is greater than for any other $I$; thus $c_\emptyset > c_I$ for any other $I$ except $I=\{1, ...., M\}$. In general, when the number of factor levels is fixed for each factor, $c_I$ is lower when the sum of the numbers of levels of the common factors in $I$ is higher ($I = \{1, ..., M\}$ excluded) and goes from positive to negative when the sum increases.

When the number of factors $M$ increases, $c_I \to 0$ for all $I \neq \{1, ..., M\}$ and the matrix $Cor$ approaches the identity matrix, that is, resembles the asymptotic correlation matrix for stratified randomization.

When $M$ remains fixed and the number of levels for factors 2 through $M$ remains fixed but $n_1$ increases, then

$$\text{for } I: 1 \in I, \ c_I \to -\frac{1}{\prod_{i=2}^{M} n_i}, \text{ while } c_I \to 0 \text{ for } I: 1 \notin I. \tag{36}$$

One practical application of (36) is an example of a multicenter study with large number of centers. When center is one of the factors in the Pocock and Simon covariate-adaptive procedure, (36) shows that under equal strata prevalence the imbalances within strata that belong to different centers are almost uncorrelated.

When the numbers of levels for at least two factors infinitely increase, from (35), $c_I \to 0$ for all $I \neq \{1, ..., M\}$.

For a common situation when all factors have 2 levels,

$$c_I = \frac{M - 1 - 2 \times \#I}{2^M - M - 1}$$

From (35), the maximum eigenvalue of the asymptotic correlation matrix $Cor$ can be derived based on the following theorem.

Theorem 2. Consider a random vector $Y$ on the set of strata $z=(i_1, ...., i_M)$, $i_j=1, ..., n_j$, such that its correlation matrix $C$ is invariant with respect to any permutation of factor levels $i_j=1, ..., n_j$ for any $j$. Then the complete set of eigenvalues of the correlation matrix $C$ is

$\{\lambda_J : J \subset \{1, ..., M\}\}$, where

$$\lambda_J = \sum_{I \subset \{1,...M\}} (-1)^{\#J^C \cap I^C} \times \tilde{c}_I \times \prod_{j \in J \cap I^C} (n_j - 1), \tag{37}$$

where $\tilde{c}_I$ is the correlation coefficient that corresponds to the two strata that have the common levels of the factors $j \in I$.

The eigenvalue $\lambda_J$ has a multiplicity

$$m_J = \prod_{j \in J^C} (n_j - 1),$$

where $m_{\{1,..,M\}} = 1$.

The proof is based on Lemma 2 described in Appendix 2; the proof is provided in Appendix 2.

Lemma 1. When $C$ is the asymptotic correlation matric $Cor$, that is $\tilde{c}_I = c_I$ in formula (35), the maximum eigenvalue of the correlation matrix $Cor$ is

$$\lambda_{max} = \frac{\prod_{i=1}^{M} n_i}{\prod_{i=1}^{M} n_i - \sum_{i=1}^{M} n_i + M - 1} \tag{38}$$

Moreover, $\lambda_J = \lambda_{max}$ for all $J$ such that $\#J < M - 1$; otherwise $\lambda_J = 0$.

Note that when the number of factors $M$ increases or the number of levels for at least one of the factors increases, $\lambda_{max} \to 1$.

From (28), under the equal prevalence condition the maximum eigenvalue of the asymptotic covariance matrix $Cov$ for the Pocock and Simon covariate-adaptive allocation is $MEV = \sigma_z \lambda_{max}$, where $\lambda_{max}$ is described by (38).

We obtained $\sigma_z^2$ through simulations under the equal prevalence condition, a bias of 0.9, in examples that included up to 7 factors and established that $\sigma_z^2 \lambda_{max} < 1$ (Appendix 1, Tables A2, A3, and A4 for 2, 3, and 4-7 factors, respectively). Generally, we observe that for large values of the bias that are used in practice, $p \geq 0.8$, $\sigma_z^2 \lambda_{max} < 1$, $\sigma_z^2 \lambda_{max}$ is the lowest for 2 factors with 2 levels each and $\sigma_z^2 \lambda_{max} \to 1$ from below as the number of strata increases. In the simulation results reported in the paper, $\sigma_z^2 \lambda_{max}$ is practically equal to 1 in most cases except for 2 factors and 3 factors with 2 levels each. We also noticed that for bias values below 0.8 ($p$=0.6, 0.67, 0.7, 0.75), $\sigma_z^2 \lambda_{max}$ is slightly above 1 for minimization with 2 factors that have 2 levels each: it is about 1.07 for $p$=0.6 and gets lower as the bias increases. For these low values of bias, $\sigma_z^2 \lambda_{max} \to 1$ as the number of factors and/or factor levels increases (simulations not shown). Although minimization is typically used with large number of factors or factor levels where even with lower values of bias $p$, $\sigma_z^2 \lambda_{max}$ is close to 1, this observation makes us recommend using the bias $p \geq 0.8$ with minimization to ensure that the maximum eigenvalue of the asymptotic covariance matrix $Cov$ is below 1 and the minimization procedure satisfies Condition 4.

Thus, for equal strata prevalence and large values of the bias ($p \geq 0.8$), it is established – partly theoretically, partly through simulations – that the maximum eigenvalue of the asymptotic covariance matrix $Cov$ for the Pocock and Simon covariate-adaptive allocation is below 1. For unequal strata prevalence, we were not able to derive the asymptotic covariance matrix and its maximum eigenvalue theoretically, however, the simulations support the maximum eigenvalue of $Cov$ being below 1 with unequal strata prevalence. Table A5 in Appendix 1 provides the maximum eigenvalue of $Cov$ for some examples of the Pocock and Simon covariate-adaptive allocation with $p$=0.9 and two independent factors one of which has unequal prevalence of its levels. In one example of unequal prevalence of the strata in Appendix 1, Table A5 the estimate of $\sigma_z^2 \lambda_{max}$ was slightly above 1 (1.00873787) and within the variability expected in estimating $\sigma_z^2$ through simulations. We also observed that the maximum eigenvalue of the asymptotic covariance matrix $Cov$ is below 1 in examples when the factors were not independent (results not provided).

Therefore, for the Pocock and Simon covariate-adaptive allocation with equal or unequal prevalence of the strata, for every vector $G$, $\gamma = G'CovG/G'G \leq 1$ for bias $p \geq 0.8$ and thus $T_S$ and $T_L$ are conservative unless $\gamma = 1$ in which case the tests are valid. However, in contrast to the stratified randomization with $Cov = vI$, $\gamma$ and thus the degree to which the Type I error with these tests is reduced compared to the nominal level depends on vector $G$. For values of bias smaller than $p=0.8$, for example, $p=0.67$ that is used in the simulations of Ye and Shao [22], the fact that $\sigma_z^2 \lambda_{max}$ is slightly above 1 does not mean that for the vector $G$ observed in a particular study, $\gamma$ is also higher than 1. If the vector $G$ has a non-negligible component orthogonal to the space spanned by the eigenvectors that correspond to the eigenvalues equal to $\lambda_{max}$, $\gamma$ is likely to be below 1 which would lead to $T_S$ and $T_L$ being conservative – as was the case in the simulations in Ye and Shao [22].

Section 4. Robust tests and simulations

Ye and Shao [22] proposed a score test robust under the conditions of Theorem 1 and randomization procedures with either condition 2 or 3 by replacing the denominator of the score test statistic (9) with a consistent estimator of the numerator's asymptotic variance. This variance estimator is a stratified version of the robust variance in Lin and Wei [39]. To construct a score test robust against model misspecification under covariate-adaptive randomization that satisfies condition 4, we generalize the test proposed by Ye and Shao [22] by replacing the denominator of the score test statistic (9) with a consistent estimator of the numerator's asymptotic variance, $\psi + G'CovG$. This follows from Theorem 1 in this paper.

As in eq. (19) of Ye and Shao [22], let

$$\hat{O}_i = \delta_i \left\{ I_i - \frac{S_N^{(1)}(\hat{\beta}_0, X_i)}{S_N^{(0)}(\hat{\beta}_0, X_i)} \right\} - \sum_{j=1}^{N} \frac{\delta_j Y_i(X_j) \exp(\hat{\beta}_0' W_i)}{N S_N^{(0)}(\hat{\beta}_0, X_j)} \left\{ I_i - \frac{S_N^{(1)}(\hat{\beta}_0, X_j)}{S_N^{(0)}(\hat{\beta}_0, X_j)} \right\},$$

where $I_i = 0, 1$ is the treatment indicator for subject $i$, $S_N^{(r)}(\hat{\beta}_0, X_i) = N^{-1} \sum_{k=1}^{N} Y_k(X_i) \exp(\hat{\beta}_0' W_k) I_k^r$, $r = 0, 1$, and $Y_k(X_i) = I_k I(X_{k1} \geq X_i) + (1 - I_k) I(X_{k0} \geq X_i)$, $I(\cdot)$ is the indicator function. The denominator of (9) can then be expressed as (Lin and Wei [23], Ye and Shao [22]):

$$\hat{B}(0, \hat{\beta}_0) = \frac{1}{N} \sum_{j=1}^{N} \hat{O}_i^2.$$

Ye and Shao [22] have shown that when $H_0$, condition 1 and assumptions 1 to 3 are true, $\hat{B}(0, \hat{\beta}_0)$ converges in probability to $\psi + \phi$ (see proof of eq. (14) in Ye and Shao [22]). Thus, if $\phi > G'CovG$, $\hat{B}(0, \hat{\beta}_0)$ overestimates the variance of the numerator $N^{-1/2} U_\theta(0, \hat{\beta}_0)$, which is possible with both stratified randomization and the Pocock and Simon minimization and model misspecification.

Let $\hat{V}_{zj}$ be the sample variance of $\hat{O}_i$ within stratum $z$ and treatment $j$, $j = 0, 1$ and let $\hat{E}_z$ be the sample mean of $\hat{O}_i$ within stratum $z$. If the covariate-adaptive randomization satisfies condition 3, meaning that the normalized treatment imbalances within the strata formed by $Z$ are asymptotically normal and uncorrelated across the strata, Ye and Shao [22] put forward the following consistent estimator of the asymptotic variance for the score test numerator, $\psi + \nu\phi$:

$$\hat{B}_{RS}(0, \hat{\beta}_0) = \frac{1}{N}\sum_z N_z \left(\frac{1}{2}\hat{V}_{z1} + \frac{1}{2}\hat{V}_{z0} + \nu \hat{E}_z^2\right),$$

where $\hat{V}_{zj}$ and $\hat{E}_z$ are consistent estimators for $\text{Var}(O_{ij}|Z_i = z)$ and $E(O_{ij}|Z_i = z)$, respectively, and $\nu$ is the known asymptotic variance of the imbalance in treatment assignments within a stratum. From Ye and Shao [22], $N^{-1}\sum_z N_z \left(\frac{1}{2}\hat{V}_{z1} + \frac{1}{2}\hat{V}_{z0}\right)$ is a consistent estimator for $\psi$ and $N^{-1}\sum_z N_z \nu \hat{E}_z^2$ is a consistent estimator for $\nu\phi$. Thus, $\hat{B}(0, \hat{\beta}_0)$, the variance in the Lin and Wei [23] robust score test, converges to a quantity where the normalized imbalance variance in a stratum must be $\nu = 1$, which would be the case if complete randomization were to be used within a stratum.

If the normalized treatment imbalances are correlated across the strata and condition 4 is satisfied, we appropriately modify this estimator to be:

$$\hat{B}_{RS}(0, \hat{\beta}_0) = \frac{1}{N}\sum_z N_z \left(\frac{1}{2}\hat{V}_{z1} + \frac{1}{2}\hat{V}_{z0}\right) + \widehat{G'CovG},$$

with $\widehat{G'CovG}$ being a consistent estimator for $G'CovG$ and $\hat{G} = (\sqrt{N_z/N}\hat{E}_z, z = 1, \ldots, M_z)'$. The difference between the two $\hat{B}_{RS}(0, \hat{\beta}_0)$ variances above is a more general estimator for the asymptotic variance of $U_2$, $G'CovG$.

Two cases where $\widehat{G'CovG}$ is tractable are as follows. If $Cov = \nu I$, where $I$ is an identity matrix and $\nu$ is known, $\widehat{G'CovG} = \frac{1}{N}\sum_z N_z \nu \hat{E}_z^2$, and the two variance estimators are the same. In this case, the consistency of $\widehat{G'CovG}$ relies only on $\hat{E}_z$ being a consistent estimator for $E(O_{ij}|Z_i = z)$, already shown in Ye and Shao [22]. In the setting with equal prevalence of all strata and the Pocock and Simon covariate-adaptive allocation, $Cov = \sigma_z^2 Cor$ with the correlation matrix $Cor$ (35). Estimating $\sigma_z^2$ via Monte Carlo and computing $\hat{E}_z$ as in Ye and Shao [22], $\widehat{G'CovG} = \hat{G}'\hat{\sigma}_z^2 Cor\hat{G}$ is a consistent estimator for $G'CovG$ by the Continuous Mapping Theorem since both $\hat{\sigma}_z^2$ and $\hat{E}_z$ are consistent for $\sigma_z^2$ and $E(O_{ij}|Z_i = z)$, respectively.

The robust score test is:

$$T_{RS} = N^{-1/2}U_\theta(0, \hat{\beta}_0)/\hat{B}_{RS}(0, \hat{\beta}_0)^{1/2}.$$

The equivalent robust log-rank test, $T_{RL}$, is constructed similarly by setting $\exp(\hat{\beta}_0' W_i) = 1$ in $\hat{O}_i$ to form $\tilde{O}_i$ (Ye and Shao, [22]), then computing $\hat{E}_z$, $\hat{V}_{zj}$, $\hat{G}$ in base of $\tilde{O}_i$ to obtain the variance:

$$\hat{B}_{RL}(0) = \frac{1}{N}\sum_z N_z \left(\frac{1}{2}\hat{V}_{z1} + \frac{1}{2}\hat{V}_{z0}\right) + \widehat{G'CovG},$$

$$T_{RL} = N^{-1/2}U/\hat{B}_{RL}(0)^{1/2},$$

where $U$ is the numerator of the log-rank test.

Robustness to model misspecification is also needed with the stratified log-rank test that does not adjust for all strata used in randomization when those strata are in fact prognostic. Such a partially stratified log-rank test may be encountered when some covariates used during randomization are omitted from the analysis model, e.g., sites, or when several sparse strata are combined in the analysis model. In this case, the numerator of the partially stratified log-rank statistic can be thought of as a sum of individual log-rank test numerators over the set of strata used in the analysis, some of which are a combination of the strata used in randomization. Similarly to the unstratified log-rank test, the numerators that correspond to the pooled strata have their variance overestimated. Hence, additively these numerators produce a misspecified log-rank test that is also not valid but conservative [22]. Since the partially stratified log-rank test is conservative, we also consider the following robust partially stratified log-rank test:

$$T_{RPL} = N^{-1/2} \sum_{\mathbf{z}'} U_{\mathbf{z}'} /\hat{B}_{RPL}(0)^{1/2},$$

where $\{\mathbf{z}'\}$ is a set of strata used in the analysis and $U_{\mathbf{z}'}$ is the numerator of the log-rank statistic in stratum $\mathbf{z}'$. To obtain the robust variance $\hat{B}_{RPL}(0)$, $\tilde{O}_{i\mathbf{z}'}$ is computed within stratum $\mathbf{z}'$ separately. The robust variance estimator $\hat{B}_{RPL}(0)$ then has the same form as $\hat{B}_{RL}(0)$ with $\hat{E}_z$, $\hat{V}_{zj}$ and $\hat{G}$ being estimated in base of $\tilde{O}_{i\mathbf{z}'}$ rather than $\tilde{O}_i$. Of note is that while the sample sizes for computing $\hat{E}_z$ and $\hat{V}_{zj}$ are the randomization strata sample sizes, the sample sizes for computing $\tilde{O}_{i\mathbf{z}'}$ are those from the analysis strata, $\mathbf{z}'$.

The asymptotic covariance matrix $Cov$ that corresponds to the population sampled in a particular clinical trial can be estimated in different ways and will be a subject of an upcoming publication. Below we will examine the properties of the proposed robust statistics following the Pocock and Simon covariate-adaptive allocation in the setting with equal prevalence of all strata, where $\widehat{Cov} = \hat{\sigma}_z^2 Cor$, with the correlation matrix $Cor$ described by (35) and the variance $\hat{\sigma}_z^2$ estimated via Monte Carlo (Appendix Tables A2-A4).

Two simulation studies are conducted. The first simulation estimates the Type 1 error and power of $T_{RS}$, $T_S$, $T_L$, $T_{RL}$ and the stratified log-rank test correctly accounting for all strata at randomization, $T_{SL}$. These tests are applied when the working model is correct and when it is misspecified with randomization following the Pocock and Simon procedure balancing across two independent and equally distributed binary factors $\mathbf{Z}_i = (Z_{i1}, Z_{i2})$. To carry out these simulations, we consider two cases:

1) Case 1: Neither of the two factors is prognostic, and the survival time for the $i$-th subject (in months) is exponentially distributed with the hazard function:

$$\ln(\lambda(t, I_i)) = \ln(0.0625) + \theta I_i.$$

2) Case 2: Both factors are highly prognostic, and the survival time for the $i$-th subject (in months) is exponentially distributed with the hazard function that follows the proportional hazard model:

$$\ln(\lambda(t, \mathbf{Z}_i, I_i)) = \ln(0.0625) + \ln(10)\, I(Z_{i1} = 2) + \ln(5)\, I(Z_{i2} = 2) + \theta I_i,$$

with $I(\cdot)$ being the indicator function. For each case, 5000 data sets with $N = 600$ were generated under the alternative hypothesis with $\theta = \ln(0.7)$ (which corresponds to the hazard ratio of 0.7 for treatment 1 vs treatment 0) and under the null hypothesis with $\theta = 0$. The sequence of independent vectors of two binary covariates was generated for the study subjects; the two covariates were independent with equally distributed levels. The Pocock and Simon covariate-adaptive randomization with bias parameter 0.9 that balances on the two factors was used to generate the treatment assignments for the study subjects. The survival time was generated following the proportional hazard model above. The intercept in the hazard function was chosen to provide a median survival of 11.1 months in the control group for the lowest risk stratum. Uniform enrollment over 29 months was simulated and the total study follow-up was set at 36 months. The study design and the simulated data represent an example of a hypothetical oncology study with progression-free survival as a primary endpoint. The nonadministrative exponentially distributed censoring time with hazard 0.01 month$^{-1}$ was applied; the overall censoring rate including the administrative censoring was approximately 14%.

Thus, the error rates of all tests are evaluated under the same simulated 5000 trials. For the robust tests $T_{RS}$ and $T_{RL}$, $\hat{\sigma}_z = 0.24$ (see Table 2A) and $Cor$ was given in (34):

$$Cor = \begin{bmatrix} 1 & -1 & -1 & 1 \\ -1 & 1 & 1 & -1 \\ -1 & 1 & 1 & -1 \\ 1 & -1 & -1 & 1 \end{bmatrix}.$$

The second simulation sets up a hypothetical trial using the Pocock and Simon procedure to balance across four independent and equally distributed binary factors $\mathbf{Z}_i = (Z_{i1}, Z_{i2}, Z_{i3}, Z_{i4})$, also using a bias parameter of 0.9. The survival time for the $i$-th subject (in months) is exponentially distributed with the hazard function:

$$\ln(\lambda(t, \mathbf{Z}_i, I_i)) = \ln(0.015) + \ln(3)\, I(Z_{i1} = 2) + \ln(2)\, I(Z_{i2} = 2) + \ln(2)\, I(Z_{i3} = 2) + \ln(2)\, I(Z_{i4} = 2) + \theta I_i.$$

The trial is analyzed under model misspecification where the last two factors are omitted from the analysis model using $T_S$, $T_L$, $T_{PL}$, and the robust tests $T_{RS}$ and $T_{RPL}$. Ten thousand data sets with $N=1000$ were generated under the alternative hypothesis with $\theta = \ln(0.78)$ and under the null hypothesis with $\theta = 0$ similar to the previous examples. The sequence of independent vectors of four binary covariates was generated for the study subjects; the two covariates were independent with equally distributed levels. The Pocock and Simon covariate-adaptive randomization with bias parameter 0.9 that balances on the four factors was used to generate the treatment assignments for the study subjects. The survival time was generated following the proportional hazard model above. The intercept in the hazard function was chosen to provide a median survival of 46.2 months in the control group for the lowest risk stratum. The uniform enrollment over 30 months was simulated; the total study follow-up was set at 50 months. The administrative censoring due to the data cutoff was the only censoring applied to the data. Here the study design and the simulated data represent an example of a hypothetical oncology study with overall survival as a primary endpoint.

For the robust tests $T_{RS}$ and $T_{RPL}$, $\hat{\sigma}_z = 0.68$ (see Table 4A) and $Cor$ is a 16 by 16 correlation matrix given in (35). It is of interest to evaluate the error rates under the model misspecification above.

The simulated Type 1 error and power for $T_L$, $T_{RL}$ and $T_{SL}$ in the first simulation are provided in Table 1. In this table $\widehat{\widetilde{G}}'\widehat{Cov}\widehat{\widetilde{G}}/\widehat{G}'\widehat{G}$ is an estimate of $\tilde{\gamma} = \widetilde{G}'Cov\widetilde{G}/\tilde{\phi}$ which was introduced in (22) as the ratio of the asymptotic variances of $U_2$ when accounting for minimization to the asymptotic variance of $U_2$ under complete randomization. By Corollary 2.1, if $\tilde{\phi} = 0$ or $\tilde{\gamma} = 1$, the test is valid; otherwise the test is conservative. The lower the $\tilde{\gamma}$, the more conservative the test is. The quantity $\widehat{\widetilde{G}}'\widehat{\widetilde{G}}/N^{-1}\sum_z N_z \left(\frac{1}{2}\hat{V}_{z1} + \frac{1}{2}\hat{V}_{z0}\right)$ in Table 1 assesses the relative magnitude of $\tilde{\phi}$ compared to $\tilde{\psi}$. The estimate for $N\tilde{\psi}$ is given by the column with $\sum_z N_z \left(\frac{1}{2}\hat{V}_{z1} + \frac{1}{2}\hat{V}_{z0}\right)$. The table also reports $N$ times the variance in the denominator of the log-rank test, $T_L$, and the robust variance for the log-rank test, $\hat{B}_{RL}(0)$. The first two rows of this table contain the results for Case 1, where no adjustment for covariates is necessary and thus the working model is correct. As stated in Ye and Shao [22], if the working model is correctly specified $\tilde{\phi} = var\{E(\tilde{O}_{ij}|Z_i)\} = \widetilde{G}'\widetilde{G} = 0$, which is observed in this simulation by noting that $\widehat{\widetilde{G}}'\widehat{\widetilde{G}}/N^{-1}\sum_z N_z \left(\frac{1}{2}\hat{V}_{z1} + \frac{1}{2}\hat{V}_{z0}\right) \approx 0$ in the first two rows of Table 1. In the scenario of Case 1, $\widetilde{G}'Cov\widetilde{G} = \tilde{\phi} = 0$ and so $\tilde{\gamma}$ is not defined, although in the presence of sampling variability one is still able to compute the ratio $\widehat{\widetilde{G}}'\widehat{Cov}\widehat{\widetilde{G}}/\widehat{G}'\widehat{G}$ (see the respective column). As prescribed by Corollary 2.1, with $\tilde{\phi} = 0$ the unstratified log-rank test preserves the nominal Type 1 error. The stratified log-rank test, which adjusts for all 4 strata used in the randomization, also preserves the Type 1 error as expected, although the covariate adjustment is not necessary in Case 1. In this case the robust variance for the log-rank test, $\hat{B}_{RL}(0)$ and the variance in the denominator of the log-rank test are very close in medians – and both are very close to the median of the estimate $N^{-1}\sum_z N_z \left(\frac{1}{2}\hat{V}_{z1} + \frac{1}{2}\hat{V}_{z0}\right)$ of $\tilde{\psi}$.

In Case 2, one can see from Table 1 that $\widehat{\widetilde{G}}'\widehat{\widetilde{G}}>0$ and $\widehat{\widetilde{G}}'\widehat{Cov}\widehat{\widetilde{G}}/\widehat{\widetilde{G}}'\widehat{\widetilde{G}}$ is very small which means that the asymptotic variance of $U_2$ is much overestimated by $\widehat{\widetilde{G}}'\widehat{\widetilde{G}}$. Since the table shows that $\widehat{\widetilde{G}}'\widehat{\widetilde{G}}$ is comparable to the estimated asymptotic variance of $U_1$ (the median of the ratio presented in the 3rd column is approximately 1.24), such gross overestimation leads to the substantial overestimation of the variance of the numerator of the log-rank test. This is turn leads to the log-rank test being very conservative (Type I error is 0.0032) and considerably less powerful than both $T_{SL}$ and $T_{RL}$. The robust log-rank test, $T_{RL}$, intends to correct the overestimation of the variance of the numerator – as can be concluded from the table, the median of the robust estimate of the variance (column $N\hat{B}_{RL}(0)$) is less than the half of the estimate of the variance of the numerator under the complete randomization, which is the column with $N$ times the variance in the denominator of the log-rank test $T_L$. The Type I error with the robust test is very close to the nominal level. However, while $T_{RL}$ preserves the Type 1 error by correcting the variance, it is less powerful than $T_{SL}$, which is consistent with the conclusion in Ye and Shao [22] that the covariate information incorporated via randomization cannot fully recover the efficiency loss due to mismodelling.

The simulated Type 1 error and power for $T_S$ and $T_{RS}$ for the first simulation when the model is correctly specified are provided in Table 2. The working models have the same functional form as the true hazard function. Hence, $\widehat{G}'\widehat{G} \approx 0$, the variance of the numerator of the score test is estimated correctly and $T_S$ preserves the Type 1 error as expected. Also, when the working model is correct and all prognostic covariates are discrete, $T_S$, $T_{RS}$ and $T_{SL}$ all have similar power (power row for Case 1 in Table 1 and power rows for Cases 1 and 2 in Table 2).

In Table 3 we report the results when fitting a misspecified model for Case 2 of the first simulation where the second covariate is omitted from the model. Specifically, the working model applied to Case 2 is:

$$\log(\lambda(t, \boldsymbol{W}_i, I_i)) = \hat{\beta}_0 + \hat{\beta}_1 I(Z_{i1} = 2) + \hat{\theta} I_i.$$

Since the model is incorrectly specified, there is no guarantee $T_S$ will preserve the Type 1 error which is the outcome in this case. Hence, with $\widehat{G}'\widehat{G} > 0$, the ratio $\widehat{G}'\widehat{Cov}\widehat{G}/\widehat{G}'\widehat{G}$ is very small and so $T_S$ is conservative and less powerful than $T_{RS}$ which preserves the Type 1 error (last column in Table 3). In this example $\widehat{G}'\widehat{G}$ is smaller relative to the estimated asymptotic variance of $U_1$ (the median of the ratio presented in the 3rd column is approximately 0.42) and thus overestimation of the asymptotic variance of $U_2$ leads to a less dramatic overestimation of the variance of the numerator than in the log-rank test example (Table 2). Also consistent with the conclusions in Ye and Shao [22], when the model is misspecified $T_{RS}$ is more powerful than $T_{RL}$ since controlling for at least some covariates improves efficiency (compare power of these tests in the last rows in Tables 1 and 3).

Table 4 reports the simulated error rates for the second simulation using a hypothetical trial with four prognostic covariates used in minimization and model misspecification through the partially stratified logrank test used to analyze the data. Here we note that similarly to the unstratified log-rank test, $T_L$, the partially stratified log-rank test, $T_{PL}$, is still conservative but more powerful than $T_L$. In fact, it behaves similarly to the score test from the same misspecified model. The respective robust score test and robust partially stratified log-rank test recover the Type 1 error through variance correction. We conclude that all simulation results in this section are consistent with the theoretical findings in this paper and with the simulation results in Ye and Shao [22] about the performance of the log-rank and score tests following the Pocock and Simon covariate-adaptive procedure.

Discussion

Unknown theoretical properties of the log-rank test and score test following minimization generated a distrust in some regulatory agencies and requests for the re-randomization test as a primary or a sensitivity analysis. In some circumstances, the re-randomization test might indeed have better performance than the stratified log-rank test that does not incorporate all baseline covariates used in minimization as stratification factors. However, the re-randomization test is computationally demanding, especially at an interim analysis when a $p$-value required for statistical significance is very small and a very large number of re-randomizations is required to estimate it with sufficient precision. In some instances re-randomization test is not clearly defined (e.g., non-inferiority tests) or requires conditioning (as with a comparison of two arms in a study with randomization to multiple arms). Thus, better understanding of the theoretical properties of the log-rank test and the score test following minimization is needed.

We advanced this understanding by expanding the methodology developed by Ye and Shao [22] for stratified randomization to covariate-adaptive randomization where the imbalances in different strata are correlated. We showed that for a class of randomization procedures that satisfy condition 4, the score test and the log-rank test are conservative under model misspecification. We showed, in part theoretically, in part through simulations, that the Pocock and Simon covariate-adaptive procedure with bias $p \geq 0.8$ belongs to this class and thus, leads to conservative log-rank and score tests under model misspecification. Although restricted to bias $p \geq 0.8$, this conclusion is very important because bias of 0.8-0.95 is typically used in practice to provide a good marginal balance in treatment assignments at any point of time. Having $p \geq 0.8$ is a sufficient condition for these tests to be conservative through Condition 4, but not a necessary condition. As was explained in most situations bias $p < 0.8$ will also lead to conservative score test and the log-rank test, however, we would recommend using bias $p \geq 0.8$ with minimization when these tests are planned. It could be further explored why lower bias leads to a slightly higher asymptotic variance of within-stratum normalized imbalances in case of equal prevalence of the strata as our simulations show.

A subsequent problem is constructing hypothesis tests that are robust against both model misspecification and randomization following minimization. Although the sandwich variance ([40], p. 160) derived by Lin and Wei [23] for the score test in the Cox model was originally constructed to be robust against model misspecification, Ye and Shao [22] show that it does not account for the correct within stratum variance of the treatment imbalance induced by stratified randomization procedures, such as the stratified permuted blocks. This impacts the hypothesis test when the model is misspecified by reducing its Type 1 error. The solution proposed by Ye and Shao [22] is a stratified version of the Lin and Wei [39] variance estimator. We generalized this robust stratified covariance estimator to accommodate randomization procedures satisfying condition 4, such as minimization, where the treatment allocations and the treatment imbalances are correlated both within and across all strata formed by the covariates used during randomization. Our simulations applying this robust covariance estimator confirm that when the model is misspecified under minimization the Type 1 error loss is recovered. Prior knowledge of the asymptotic covariance matrix of the normalized within stratum imbalances is necessary to apply this robust variance estimator.

While the substantial progress was made in understanding the properties of the score test and the log-rank tests following the Pocock and Simon covariate-adaptive randomization, open questions remain. When the strata have equal prevalence, we show that under minimization the asymptotic covariance matrix is the product of the common asymptotic variance of the normalized within stratum imbalances and the asymptotic correlation matrix of the normalized within stratum imbalances described by (35). We derived theoretically the asymptotic correlation matrix for the case of $M=2$ factors with any number of levels; however, for $M>2$ factors we supported the formula (35) by simulation results. Obtaining the theoretical proof of (35) for more than 2 factors would be a worthy pursuit. We derived theoretically the maximum eigenvalue of the correlation matrix (35) in Theorem 2 and Lemma 1. We used the estimates of the variance of the normalized within-stratum imbalances obtained through simulations to show that the product of the variance and the maximum eigenvalue – that is, the maximum eigenvalue of the asymptotic covariance function – is below 1 with bias $p \geq 0.8$. A theoretical derivation of the variance of the normalized within-stratum imbalances in the equal prevalence setting would be valuable as well as the study of the impact of the bias $p$ on this variance.

In the general case of unequally prevalent strata, we showed through simulations in different settings that, similar to the equal prevalence case, with $p \geq 0.8$ the maximum eigenvalue of the asymptotic covariance matrix is below 1. However, the asymptotic covariance matrix of the normalized within-stratum imbalances and a theoretical demonstration that its maximum eigenvalue is less than 1 in this setting remain open questions. Resolving these open problems would constitute a fully theoretical proof that the Pocock and Simon covariate-adaptive procedure satisfies condition 4 with $p \geq 0.8$.

Acknowledgements

This research did not receive any specific grant from funding agencies in the public, commercial, or not-for-profit sectors.

**Table 1.** Simulated Type 1 error and power for the one-sided $\alpha = 0.025$ log-rank Test ($T_L$), the log-rank test with robust variance ($T_{RL}$) and the stratified log-rank Test ($T_{SL}$) in 5,000 simulation runs, N=600.

| True Hazard Scenario | Unstratified log-rank test ($T_L$) | | | | | | $T_{RL}$ Type 1 error or power | $T_{SL}$ Type 1 error or power |
|---|---|---|---|---|---|---|---|---|
| | Type 1 error or power | $\widehat{G}'\widehat{Cov}\widehat{G}/\widehat{G}'\widehat{G}$ median (max, min) | $\widehat{G}'\widehat{G}/N^{-1}\sum_z N_z \left(\frac{1}{2}\widehat{V}_{z1} + \frac{1}{2}\widehat{V}_{z0}\right)$ median (max, min) | $\sum_z N_z \left(\frac{1}{2}\widehat{V}_{z1} + \frac{1}{2}\widehat{V}_{z0}\right)$ median (max, min) | $N$ times the variance in the denominator of $T_L$ median (max, min) | $N\widehat{B}_{RL}(0)$ median (max, min) | | |
| Case 1: Type I error | 0.0244 | 0.2386 (<0.0001, 0.9588) | 0.0040 (<0.0001, 0.0305) | 96.1297 (86.0497, 106.2683) | 96.4366 (86.6409, 105.9820) | 96.2631 (86.0587, 106.2899) | 0.0250 | 0.0248 |
| Case 1: Power | 0.9118 | 0.2368 (<0.0001, 0.9558) | 0.0042 (<0.0001, 0.0327) | 84.9293 (72.4039, 95.9987) | 87.0392 (75.1757, 97.0264) | 85.0455 (72.4318, 96.2442) | 0.9138 | 0.9078 |
| Case 2: Type I error | 0.0032 | 0.0628 (0.0188, 0.1259) | 1.2379 (0.8507, 1.8268) | 59.7181 (44.8545, 73.9773) | 133.9290 (127.7950, 140.8600) | 64.4555 (49.3734, 77.6209) | 0.0242 | 0.0270 |
| Case 2: Power | 0.5972 | 0.0584 (0.0164, 0.1244) | 1.2347 (0.8190, 1.9142) | 57.1480 (42.9162, 73.0765) | 130.1270 (121.9360, 137.5330) | 61.3901 (46.9138, 76.6120) | 0.8716 | 0.9802 |

**Table 2.** Correct Model: Simulated Type 1 error and power for the one-sided $\alpha = 0.025$ score test ($T_S$) and the robust score test ($T_{RS}$) in 5,000 simulation runs, N=600.

| True Hazard Scenario | Score test ($T_S$) | | | | | | $T_{RS}$ Type I error or power |
|---|---|---|---|---|---|---|---|
| | Type 1 error or power | $\widehat{G}'\widehat{Cov}\widehat{G}/\widehat{G}'\widehat{G}$ median (max, min) | $\widehat{G}'\widehat{G}/N^{-1}\sum_z N_z\left(\frac{1}{2}\widehat{V}_{z1}+\frac{1}{2}\widehat{V}_{z0}\right)$ median (max, min) | $\sum_z N_z\left(\frac{1}{2}\widehat{V}_{z1}+\frac{1}{2}\widehat{V}_{z0}\right)$ median (max, min) | $N\widehat{B}(0,\hat{\beta}_0)$ median (max, min) | $N\widehat{B}_{RS}(0,\hat{\beta}_0)$ median (max, min) | |
| Case 1: Type I error | 0.0244 | 0.2368 (<0.0001, 0.9503) | 0.0040 (<0.0001, 0.0305) | 96.1297 (86.0497, 106.2683) | 95.9182 (85.8611, 105.8810) | 96.2631 (86.0587, 106.2899) | 0.0248 |
| Case 1: Power | 0.9130 | 0.2368 (<0.0001, 0.9558) | 0.0042 (<0.0001, 0.0327) | 84.9293 (72.4039, 95.9987) | 86.2991 (75.0558, 96.2666) | 85.0455 (72.0558, 96.2442) | 0.9138 |
| Case 2: Type I error | 0.0270 | 0.9200 (<0.0001, 0.9599) | 0.0007 (<0.0001, 0.0190) | 132.0514 (111.7517, 160.5872) | 131.2738 (111.2241, 158.9674) | 132.2078 (111.9590, 160.7418) | 0.0274 |
| Case 2: Power | 0.9824 | 0.4585 (<0.0001, 0.9584) | 0.0017 (<0.0001, 0.0230) | 119.1069 (93.2304, 146.0638) | 121.8008 (101.8054, 147.8480) | 119.2419 (93.6243, 146.0961) | 0.9820 |

**Table 3.** Incorrect Model: Simulated Type 1 error and power for the one-sided $\alpha = 0.025$ score test ($T_S$) and the robust score test ($T_{RS}$) in 5,000 simulation runs, N=600.

| True Hazard Scenario | Score Test ($T_S$) | | | | | | $T_{RS}$ Type I error or power |
|---|---|---|---|---|---|---|---|
| | Type 1 error or power | $\widehat{G}'\widehat{Cov}\widehat{G}/\widehat{G}'\widehat{G}$ median (max, min) | $\widehat{G}'\widehat{G}/N^{-1}\sum_z N_z \left(\frac{1}{2}\widehat{V}_{z1} + \frac{1}{2}\widehat{V}_{z0}\right)$ median (max, min) | $\sum_z N_z \left(\frac{1}{2}\widehat{V}_{z1} + \frac{1}{2}\widehat{V}_{z0}\right)$ median (max, min) | $N\widehat{B}(0,\hat{\beta}_0)$ median (max, min) | $N\widehat{B}_{RS}(0,\hat{\beta}_0)$ median (max, min) | |
| Case 2: Type I error | 0.0096 | 0.0031 (<0.0001, 0.0524) | 0.4186 (0.2726, 0.6272) | 98.6195 (78.9304, 123.2212) | 139.0959 (123.9973, 159.4351) | 98.8438 (78.9743, 123.8836) | 0.0256 |
| Case 2: Power | 0.8858 | 0.0038 (<0.0001, 0.0623) | 0.4185 (0.2531, 0.6295) | 91.8034 (70.5158, 118.6563) | 131.9277 (114.2020, 152.6467) | 91.9967 (70.5189, 118.7273) | 0.9436 |

**Table 4:** Simulated Type 1 error and power for the one-sided $\alpha = 0.025$ score test ($T_S$), the robust score test ($T_{RS}$), log-rank test ($T_L$), partially stratified log-rank test ($T_{PL}$), and the robust partially stratified log-rank test ($T_{RPL}$) following a misspecified model using 10,000 simulation runs, N=1000.

| Testing procedure | Type 1 error rate | Power |
|---|---|---|
| $T_S$ | 0.0168 | 0.8871 |
| $T_{RS}$ | 0.0267 | 0.9141 |
| $T_L$ | 0.0093 | 0.7911 |
| $T_{PL}$ | 0.0165 | 0.8864 |
| $T_{RPL}$ | 0.0260 | 0.9138 |

Appendix 1.

Table A1. Theoretical correlations of the within-stratum imbalances in treatment assignments under equal prevalence (35) vs. correlations derived through simulations: 4-factor examples of Pocock and Simon covariate-adaptive allocation with $p=0.9$, 10,000 simulation runs and $500 \times M_s$ subjects per run.

| $n_1$ | $n_2$ | $n_3$ | $n_4$ | $\varepsilon_1$ | $\varepsilon_2$ | $\varepsilon_3$ | $\varepsilon_4$ | Theoretical correlation | Correlation derived through simulations | Absolute Difference Between Theoretical and simulated correlations |
|---|---|---|---|---|---|---|---|---|---|---|
| 2 | 2 | 2 | 2 | 0 | 0 | 0 | 0 | 0.27273 | 0.26202 | 0.01071 |
| 2 | 2 | 2 | 2 | 0 | 0 | 0 | 1 | 0.09091 | 0.09513 | 0.00422 |
| 2 | 2 | 2 | 2 | 0 | 0 | 1 | 0 | 0.09091 | 0.08669 | 0.00422 |
| 2 | 2 | 2 | 2 | 0 | 0 | 1 | 1 | -0.09091 | -0.08652 | 0.00439 |
| 2 | 2 | 2 | 2 | 0 | 1 | 0 | 0 | 0.09091 | 0.10398 | 0.01307 |
| 2 | 2 | 2 | 2 | 0 | 1 | 0 | 1 | -0.09091 | -0.09772 | 0.00681 |
| 2 | 2 | 2 | 2 | 0 | 1 | 1 | 0 | -0.09091 | -0.08703 | 0.00388 |
| 2 | 2 | 2 | 2 | 0 | 1 | 1 | 1 | -0.27273 | -0.27691 | 0.00419 |
| 2 | 2 | 2 | 2 | 1 | 0 | 0 | 0 | 0.09091 | 0.09037 | 0.00054 |
| 2 | 2 | 2 | 2 | 1 | 0 | 0 | 1 | -0.09091 | -0.08661 | 0.00430 |
| 2 | 2 | 2 | 2 | 1 | 0 | 1 | 0 | -0.09091 | -0.09116 | 0.00025 |
| 2 | 2 | 2 | 2 | 1 | 0 | 1 | 1 | -0.27273 | -0.27029 | 0.00244 |
| 2 | 2 | 2 | 2 | 1 | 1 | 0 | 0 | -0.09091 | -0.09152 | 0.00061 |
| 2 | 2 | 2 | 2 | 1 | 1 | 0 | 1 | -0.27273 | -0.27601 | 0.00328 |
| 2 | 2 | 2 | 2 | 1 | 1 | 1 | 0 | -0.27273 | -0.27370 | 0.00097 |
| 2 | 2 | 2 | 3 | 0 | 0 | 0 | 0 | 0.16667 | 0.16605 | 0.00062 |
| 2 | 2 | 2 | 3 | 0 | 0 | 0 | 1 | 0.00000 | -0.00291 | 0.00291 |
| 2 | 2 | 2 | 3 | 0 | 0 | 1 | 0 | 0.05556 | 0.05365 | 0.00191 |
| 2 | 2 | 2 | 3 | 0 | 0 | 1 | 1 | -0.11111 | -0.11785 | 0.00673 |
| 2 | 2 | 2 | 3 | 0 | 1 | 0 | 0 | 0.05556 | 0.05686 | 0.00131 |
| 2 | 2 | 2 | 3 | 0 | 1 | 0 | 1 | -0.11111 | -0.10378 | 0.00734 |
| 2 | 2 | 2 | 3 | 0 | 1 | 1 | 0 | -0.05556 | -0.05517 | 0.00039 |
| 2 | 2 | 2 | 3 | 0 | 1 | 1 | 1 | -0.22222 | -0.21845 | 0.00378 |
| 2 | 2 | 2 | 3 | 1 | 0 | 0 | 0 | 0.05556 | 0.05443 | 0.00113 |
| 2 | 2 | 2 | 3 | 1 | 0 | 0 | 1 | -0.11111 | -0.10632 | 0.00479 |
| 2 | 2 | 2 | 3 | 1 | 0 | 1 | 0 | -0.05556 | -0.04866 | 0.00689 |
| 2 | 2 | 2 | 3 | 1 | 0 | 1 | 1 | -0.22222 | -0.22404 | 0.00182 |
| 2 | 2 | 2 | 3 | 1 | 1 | 0 | 0 | -0.05556 | -0.05818 | 0.00262 |
| 2 | 2 | 2 | 3 | 1 | 1 | 0 | 1 | -0.22222 | -0.22552 | 0.00330 |
| 2 | 2 | 2 | 3 | 1 | 1 | 1 | 0 | -0.16667 | -0.16921 | 0.00255 |
| 2 | 2 | 2 | 4 | 0 | 0 | 0 | 0 | 0.12000 | 0.12181 | 0.00181 |
| 2 | 2 | 2 | 4 | 0 | 0 | 0 | 1 | -0.04000 | -0.04067 | 0.00067 |

| 2 | 2 | 2 | 4 | 0 | 0 | 1 | 0 | 0.04000 | 0.03854 | 0.00146 |
|---|---|---|---|---|---|---|---|---|---|---|
| 2 | 2 | 2 | 4 | 0 | 0 | 1 | 1 | -0.12000 | -0.12218 | 0.00218 |
| 2 | 2 | 2 | 4 | 0 | 1 | 0 | 0 | 0.04000 | 0.03742 | 0.00258 |
| 2 | 2 | 2 | 4 | 0 | 1 | 0 | 1 | -0.12000 | -0.12094 | 0.00094 |
| 2 | 2 | 2 | 4 | 0 | 1 | 1 | 0 | -0.04000 | -0.03778 | 0.00222 |
| 2 | 2 | 2 | 4 | 0 | 1 | 1 | 1 | -0.20000 | -0.19634 | 0.00366 |
| 2 | 2 | 2 | 4 | 1 | 0 | 0 | 0 | 0.04000 | 0.04081 | 0.00081 |
| 2 | 2 | 2 | 4 | 1 | 0 | 0 | 1 | -0.12000 | -0.11607 | 0.00393 |
| 2 | 2 | 2 | 4 | 1 | 0 | 1 | 0 | -0.04000 | -0.04056 | 0.00056 |
| 2 | 2 | 2 | 4 | 1 | 0 | 1 | 1 | -0.20000 | -0.20304 | 0.00304 |
| 2 | 2 | 2 | 4 | 1 | 1 | 0 | 0 | -0.04000 | -0.04099 | 0.00099 |
| 2 | 2 | 2 | 4 | 1 | 1 | 0 | 1 | -0.20000 | -0.19964 | 0.00036 |
| 2 | 2 | 2 | 4 | 1 | 1 | 1 | 0 | -0.12000 | -0.11943 | 0.00057 |
| 2 | 2 | 2 | 5 | 0 | 0 | 0 | 0 | 0.09375 | 0.09389 | 0.00014 |
| 2 | 2 | 2 | 5 | 0 | 0 | 0 | 1 | -0.06250 | -0.06292 | 0.00042 |
| 2 | 2 | 2 | 5 | 0 | 0 | 1 | 0 | 0.03125 | 0.03160 | 0.00035 |
| 2 | 2 | 2 | 5 | 0 | 0 | 1 | 1 | -0.12500 | -0.12186 | 0.00314 |
| 2 | 2 | 2 | 5 | 0 | 1 | 0 | 0 | 0.03125 | 0.03078 | 0.00047 |
| 2 | 2 | 2 | 5 | 0 | 1 | 0 | 1 | -0.12500 | -0.12324 | 0.00176 |
| 2 | 2 | 2 | 5 | 0 | 1 | 1 | 0 | -0.03125 | -0.03258 | 0.00133 |
| 2 | 2 | 2 | 5 | 0 | 1 | 1 | 1 | -0.18750 | -0.18687 | 0.00063 |
| 2 | 2 | 2 | 5 | 1 | 0 | 0 | 0 | 0.03125 | 0.02963 | 0.00162 |
| 2 | 2 | 2 | 5 | 1 | 0 | 0 | 1 | -0.12500 | -0.12726 | 0.00226 |
| 2 | 2 | 2 | 5 | 1 | 0 | 1 | 0 | -0.03125 | -0.02969 | 0.00156 |
| 2 | 2 | 2 | 5 | 1 | 0 | 1 | 1 | -0.18750 | -0.18984 | 0.00234 |
| 2 | 2 | 2 | 5 | 1 | 1 | 0 | 0 | -0.03125 | -0.02926 | 0.00199 |
| 2 | 2 | 2 | 5 | 1 | 1 | 0 | 1 | -0.18750 | -0.18685 | 0.00065 |
| 2 | 2 | 2 | 5 | 1 | 1 | 1 | 0 | -0.09375 | -0.09451 | 0.00076 |
| 2 | 2 | 2 | 6 | 0 | 0 | 0 | 0 | 0.07692 | 0.07559 | 0.00133 |
| 2 | 2 | 2 | 6 | 0 | 0 | 0 | 1 | -0.07692 | -0.08017 | 0.00324 |
| 2 | 2 | 2 | 6 | 0 | 0 | 1 | 0 | 0.02564 | 0.02468 | 0.00096 |
| 2 | 2 | 2 | 6 | 0 | 0 | 1 | 1 | -0.12821 | -0.12063 | 0.00757 |
| 2 | 2 | 2 | 6 | 0 | 1 | 0 | 0 | 0.02564 | 0.02767 | 0.00203 |
| 2 | 2 | 2 | 6 | 0 | 1 | 0 | 1 | -0.12821 | -0.12937 | 0.00117 |
| 2 | 2 | 2 | 6 | 0 | 1 | 1 | 0 | -0.02564 | -0.02647 | 0.00083 |
| 2 | 2 | 2 | 6 | 0 | 1 | 1 | 1 | -0.17949 | -0.17727 | 0.00222 |
| 2 | 2 | 2 | 6 | 1 | 0 | 0 | 0 | 0.02564 | 0.02685 | 0.00121 |
| 2 | 2 | 2 | 6 | 1 | 0 | 0 | 1 | -0.12821 | -0.12788 | 0.00033 |
| 2 | 2 | 2 | 6 | 1 | 0 | 1 | 0 | -0.02564 | -0.02531 | 0.00033 |
| 2 | 2 | 2 | 6 | 1 | 0 | 1 | 1 | -0.17949 | -0.18049 | 0.00100 |
| 2 | 2 | 2 | 6 | 1 | 1 | 0 | 0 | -0.02564 | -0.02604 | 0.00040 |
| 2 | 2 | 2 | 6 | 1 | 1 | 0 | 1 | -0.17949 | -0.18304 | 0.00356 |
| 2 | 2 | 2 | 6 | 1 | 1 | 1 | 0 | -0.07692 | -0.07709 | 0.00016 |
| 2 | 2 | 3 | 3 | 0 | 0 | 0 | 0 | 0.10345 | 0.10329 | 0.00015 |
| 2 | 2 | 3 | 3 | 0 | 0 | 0 | 1 | 0.00000 | -0.00144 | 0.00144 |
| 2 | 2 | 3 | 3 | 0 | 0 | 1 | 0 | 0.00000 | 0.00040 | 0.00040 |
| 2 | 2 | 3 | 3 | 0 | 0 | 1 | 1 | -0.10345 | -0.10572 | 0.00227 |
| 2 | 2 | 3 | 3 | 0 | 1 | 0 | 0 | 0.03448 | 0.03466 | 0.00018 |
| 2 | 2 | 3 | 3 | 0 | 1 | 0 | 1 | -0.06897 | -0.06929 | 0.00033 |
| 2 | 2 | 3 | 3 | 0 | 1 | 1 | 0 | -0.06897 | -0.06885 | 0.00012 |
| 2 | 2 | 3 | 3 | 0 | 1 | 1 | 1 | -0.17241 | -0.16785 | 0.00457 |
| 2 | 2 | 3 | 3 | 1 | 0 | 0 | 0 | 0.03448 | 0.03469 | 0.00020 |
| 2 | 2 | 3 | 3 | 1 | 0 | 0 | 1 | -0.06897 | -0.06565 | 0.00332 |
| 2 | 2 | 3 | 3 | 1 | 0 | 1 | 0 | -0.06897 | -0.06943 | 0.00047 |
| 2 | 2 | 3 | 3 | 1 | 0 | 1 | 1 | -0.17241 | -0.17407 | 0.00166 |
| 2 | 2 | 3 | 3 | 1 | 1 | 0 | 0 | -0.03448 | -0.03484 | 0.00036 |
| 2 | 2 | 3 | 3 | 1 | 1 | 0 | 1 | -0.13793 | -0.13936 | 0.00143 |
| 2 | 2 | 3 | 3 | 1 | 1 | 1 | 0 | -0.13793 | -0.13785 | 0.00008 |
| 2 | 2 | 3 | 4 | 0 | 0 | 0 | 0 | 0.07500 | 1.00000 | 1.24138 |
| 2 | 2 | 3 | 4 | 0 | 0 | 0 | 1 | -0.02500 | 0.07325 | 0.00175 |
| 2 | 2 | 3 | 4 | 0 | 0 | 1 | 0 | 0.00000 | -0.02307 | 0.00193 |
| 2 | 2 | 3 | 4 | 0 | 0 | 1 | 1 | -0.10000 | 0.00018 | 0.00018 |
| 2 | 2 | 3 | 4 | 0 | 1 | 0 | 0 | 0.02500 | -0.09964 | 0.00036 |
| 2 | 2 | 3 | 4 | 0 | 1 | 0 | 1 | -0.07500 | 0.02571 | 0.00071 |
| 2 | 2 | 3 | 4 | 0 | 1 | 1 | 0 | -0.05000 | -0.07384 | 0.00116 |
| 2 | 2 | 3 | 4 | 0 | 1 | 1 | 1 | -0.15000 | -0.05001 | 0.00001 |
| 2 | 2 | 3 | 4 | 1 | 0 | 0 | 0 | 0.02500 | -0.15092 | 0.00092 |
| 2 | 2 | 3 | 4 | 1 | 0 | 0 | 1 | -0.07500 | 0.02520 | 0.00020 |
| 2 | 2 | 3 | 4 | 1 | 0 | 1 | 0 | -0.05000 | -0.07538 | 0.00038 |
| 2 | 2 | 3 | 4 | 1 | 0 | 1 | 1 | -0.15000 | -0.04782 | 0.00218 |
| 2 | 2 | 3 | 4 | 1 | 1 | 0 | 0 | -0.02500 | -0.15134 | 0.00134 |
| 2 | 2 | 3 | 4 | 1 | 1 | 0 | 1 | -0.12500 | -0.02466 | 0.00034 |
| 2 | 2 | 3 | 4 | 1 | 1 | 1 | 0 | -0.10000 | -0.12631 | 0.00131 |
| 2 | 2 | 3 | 5 | 0 | 0 | 0 | 0 | 0.05882 | -0.10146 | 0.00146 |
| 2 | 2 | 3 | 5 | 0 | 0 | 0 | 1 | -0.03922 | 0.05847 | 0.00035 |
| 2 | 2 | 3 | 5 | 0 | 0 | 1 | 0 | 0.00000 | -0.03917 | 0.00005 |
| 2 | 2 | 3 | 5 | 0 | 0 | 1 | 1 | -0.09804 | 0.00105 | 0.00105 |
| 2 | 2 | 3 | 5 | 0 | 1 | 0 | 0 | 0.01961 | -0.09537 | 0.00267 |
| 2 | 2 | 3 | 5 | 0 | 1 | 0 | 1 | -0.07843 | 0.02007 | 0.00047 |
| 2 | 2 | 3 | 5 | 0 | 1 | 1 | 0 | -0.03922 | -0.07912 | 0.00069 |
| 2 | 2 | 3 | 5 | 0 | 1 | 1 | 1 | -0.13725 | -0.04038 | 0.00117 |
| 2 | 2 | 3 | 5 | 1 | 0 | 0 | 0 | 0.01961 | -0.13914 | 0.00188 |
| 2 | 2 | 3 | 5 | 1 | 0 | 0 | 1 | -0.07843 | 0.01870 | 0.00090 |
| 2 | 2 | 3 | 5 | 1 | 0 | 1 | 0 | -0.03922 | -0.07743 | 0.00101 |
| 2 | 2 | 3 | 5 | 1 | 0 | 1 | 1 | -0.13725 | -0.03893 | 0.00028 |
| 2 | 2 | 3 | 5 | 1 | 1 | 0 | 0 | -0.01961 | -0.13735 | 0.00010 |
| 2 | 2 | 3 | 5 | 1 | 1 | 0 | 1 | -0.11765 | -0.01886 | 0.00074 |
| 2 | 2 | 3 | 5 | 1 | 1 | 1 | 0 | -0.07843 | -0.11790 | 0.00025 |
| 2 | 2 | 3 | 6 | 0 | 0 | 0 | 0 | 0.04839 | -0.07859 | 0.00015 |
| 2 | 2 | 3 | 6 | 0 | 0 | 0 | 1 | -0.04839 | 0.04763 | 0.00075 |
| 2 | 2 | 3 | 6 | 0 | 0 | 1 | 0 | 0.00000 | -0.04845 | 0.00007 |
| 2 | 2 | 3 | 6 | 0 | 0 | 1 | 1 | -0.09677 | 0.00135 | 0.00135 |
| 2 | 2 | 3 | 6 | 0 | 1 | 0 | 0 | 0.01613 | -0.09408 | 0.00269 |
| 2 | 2 | 3 | 6 | 0 | 1 | 0 | 1 | -0.08065 | 0.01637 | 0.00025 |
| 2 | 2 | 3 | 6 | 0 | 1 | 1 | 0 | -0.03226 | -0.08018 | 0.00046 |
| 2 | 2 | 3 | 6 | 0 | 1 | 1 | 1 | -0.12903 | -0.03310 | 0.00084 |
| 2 | 2 | 3 | 6 | 1 | 0 | 0 | 0 | 0.01613 | -0.13002 | 0.00099 |
| 2 | 2 | 3 | 6 | 1 | 0 | 0 | 1 | -0.08065 | 0.01612 | 0.00001 |
| 2 | 2 | 3 | 6 | 1 | 0 | 1 | 0 | -0.03226 | -0.07890 | 0.00174 |
| 2 | 2 | 3 | 6 | 1 | 0 | 1 | 1 | -0.12903 | -0.03289 | 0.00063 |
| 2 | 2 | 3 | 6 | 1 | 1 | 0 | 0 | -0.01613 | -0.13114 | 0.00210 |
| 2 | 2 | 3 | 6 | 1 | 1 | 0 | 1 | -0.11290 | -0.01576 | 0.00037 |
| 2 | 2 | 3 | 6 | 1 | 1 | 1 | 0 | -0.06452 | -0.11439 | 0.00149 |
| 2 | 2 | 4 | 4 | 0 | 0 | 0 | 0 | 0.05455 | -0.06418 | 0.00034 |
| 2 | 2 | 4 | 4 | 0 | 0 | 0 | 1 | -0.01818 | 0.05334 | 0.00121 |
| 2 | 2 | 4 | 4 | 0 | 0 | 1 | 0 | -0.01818 | -0.01749 | 0.00070 |
| 2 | 2 | 4 | 4 | 0 | 0 | 1 | 1 | -0.09091 | -0.01815 | 0.00003 |
| 2 | 2 | 4 | 4 | 0 | 1 | 0 | 0 | 0.01818 | -0.08875 | 0.00216 |
| 2 | 2 | 4 | 4 | 0 | 1 | 0 | 1 | -0.05455 | 0.01848 | 0.00030 |
| 2 | 2 | 4 | 4 | 0 | 1 | 1 | 0 | -0.05455 | -0.05275 | 0.00179 |
| 2 | 2 | 4 | 4 | 0 | 1 | 1 | 1 | -0.12727 | -0.05456 | 0.00002 |
| 2 | 2 | 4 | 4 | 1 | 0 | 0 | 0 | 0.01818 | -0.12878 | 0.00151 |
| 2 | 2 | 4 | 4 | 1 | 0 | 0 | 1 | -0.05455 | 0.01839 | 0.00021 |
| 2 | 2 | 4 | 4 | 1 | 0 | 1 | 0 | -0.05455 | -0.05537 | 0.00083 |
| 2 | 2 | 4 | 4 | 1 | 0 | 1 | 1 | -0.12727 | -0.05270 | 0.00184 |
| 2 | 2 | 4 | 4 | 1 | 1 | 0 | 0 | -0.01818 | -0.12576 | 0.00151 |
| 2 | 2 | 4 | 4 | 1 | 1 | 0 | 1 | -0.09091 | -0.01735 | 0.00083 |
| 2 | 2 | 4 | 4 | 1 | 1 | 1 | 0 | -0.09091 | -0.09303 | 0.00212 |
| 2 | 2 | 4 | 5 | 0 | 0 | 0 | 0 | 0.04286 | 0.04253 | 0.00032 |
| 2 | 2 | 4 | 5 | 0 | 0 | 0 | 1 | -0.02857 | -0.02829 | 0.00028 |
| 2 | 2 | 4 | 5 | 0 | 0 | 1 | 0 | -0.01429 | -0.01364 | 0.00065 |
| 2 | 2 | 4 | 5 | 0 | 0 | 1 | 1 | -0.08571 | -0.08496 | 0.00076 |
| 2 | 2 | 4 | 5 | 0 | 1 | 0 | 0 | 0.01429 | 0.01435 | 0.00006 |
| 2 | 2 | 4 | 5 | 0 | 1 | 0 | 1 | -0.05714 | -0.05697 | 0.00018 |
| 2 | 2 | 4 | 5 | 0 | 1 | 1 | 0 | -0.04286 | -0.04344 | 0.00058 |
| 2 | 2 | 4 | 5 | 0 | 1 | 1 | 1 | -0.11429 | -0.11357 | 0.00072 |
| 2 | 2 | 4 | 5 | 1 | 0 | 0 | 0 | 0.01429 | 0.01420 | 0.00009 |
| 2 | 2 | 4 | 5 | 1 | 0 | 0 | 1 | -0.05714 | -0.05671 | 0.00043 |
| 2 | 2 | 4 | 5 | 1 | 0 | 1 | 0 | -0.04286 | -0.04327 | 0.00041 |
| 2 | 2 | 4 | 5 | 1 | 0 | 1 | 1 | -0.11429 | -0.11327 | 0.00101 |
| 2 | 2 | 4 | 5 | 1 | 1 | 0 | 0 | -0.01429 | -0.01381 | 0.00047 |
| 2 | 2 | 4 | 5 | 1 | 1 | 0 | 1 | -0.08571 | -0.08717 | 0.00145 |
| 2 | 2 | 4 | 5 | 1 | 1 | 1 | 0 | -0.07143 | -0.07153 | 0.00010 |
| 2 | 2 | 4 | 6 | 0 | 0 | 0 | 0 | 0.03529 | 0.03519 | 0.00011 |
| 2 | 2 | 4 | 6 | 0 | 0 | 0 | 1 | -0.03529 | -0.03554 | 0.00025 |
| 2 | 2 | 4 | 6 | 0 | 0 | 1 | 0 | -0.01176 | -0.01211 | 0.00034 |
| 2 | 2 | 4 | 6 | 0 | 0 | 1 | 1 | -0.08235 | -0.08401 | 0.00166 |
| 2 | 2 | 4 | 6 | 0 | 1 | 0 | 0 | 0.01176 | 0.01187 | 0.00011 |
| 2 | 2 | 4 | 6 | 0 | 1 | 0 | 1 | -0.05882 | -0.05856 | 0.00027 |
| 2 | 2 | 4 | 6 | 0 | 1 | 1 | 0 | -0.03529 | -0.03466 | 0.00063 |
| 2 | 2 | 4 | 6 | 0 | 1 | 1 | 1 | -0.10588 | -0.10577 | 0.00011 |
| 2 | 2 | 4 | 6 | 1 | 0 | 0 | 0 | 0.01176 | 0.01180 | 0.00003 |
| 2 | 2 | 4 | 6 | 1 | 0 | 0 | 1 | -0.05882 | -0.05813 | 0.00069 |
| 2 | 2 | 4 | 6 | 1 | 0 | 1 | 0 | -0.03529 | -0.03489 | 0.00041 |
| 2 | 2 | 4 | 6 | 1 | 0 | 1 | 1 | -0.10588 | -0.10475 | 0.00113 |
| 2 | 2 | 4 | 6 | 1 | 1 | 0 | 0 | -0.01176 | -0.01189 | 0.00012 |
| 2 | 2 | 4 | 6 | 1 | 1 | 0 | 1 | -0.08235 | -0.08265 | 0.00030 |
| 2 | 2 | 4 | 6 | 1 | 1 | 1 | 0 | -0.05882 | -0.05930 | 0.00047 |

Table A2. The simulated variance of the within-stratum imbalances in treatment assignments under equal prevalence for 2-factor examples of Pocock and Simon covariate-adaptive allocation with $p=0.9$, 10,000 simulation runs and $500 \times M_s$ subjects per run; the theoretical maximum eigenvalue of the asymptotic correlation matrix of the within-stratum imbalances, and the maximum eigenvalue of the asymptotic covariance matrix of the within-stratum imbalances.

| $n_1$ | $n_2$ | Simulated variance of the within-stratum imbalances, $\sigma_z^2$ | Theoretical maximum eigenvalue of the asymptotic correlation matrix, $\lambda_{max}$ | Maximum eigenvalue of the asymptotic covariance matrix, $\sigma_z^2 \lambda_{max}$ |
|---|---|---|---|---|
| 2 | 2 | 0.23509 | 4.00000 | 0.94035 |
| 2 | 3 | 0.32176 | 3.00000 | 0.96528 |
| 2 | 4 | 0.36708 | 2.66667 | 0.97889 |
| 2 | 5 | 0.38949 | 2.50000 | 0.97373 |
| 2 | 6 | 0.40777 | 2.40000 | 0.97866 |
| 2 | 7 | 0.41738 | 2.33333 | 0.97388 |
| 2 | 8 | 0.43064 | 2.28571 | 0.98431 |
| 3 | 3 | 0.43068 | 2.25000 | 0.96904 |
| 3 | 4 | 0.48277 | 2.00000 | 0.96554 |
| 3 | 5 | 0.51880 | 1.87500 | 0.97276 |
| 3 | 6 | 0.54057 | 1.80000 | 0.97303 |
| 3 | 7 | 0.56153 | 1.75000 | 0.98268 |
| 3 | 8 | 0.57381 | 1.71429 | 0.98367 |
| 4 | 4 | 0.54804 | 1.77778 | 0.97429 |
| 4 | 5 | 0.59126 | 1.66667 | 0.98544 |
| 4 | 6 | 0.61087 | 1.60000 | 0.97739 |
| 4 | 7 | 0.62693 | 1.55556 | 0.97522 |
| 4 | 8 | 0.64688 | 1.52381 | 0.98573 |
| 5 | 5 | 0.63050 | 1.56250 | 0.98515 |
| 5 | 6 | 0.65532 | 1.50000 | 0.98297 |
| 5 | 7 | 0.67480 | 1.45833 | 0.98409 |
| 5 | 8 | 0.68902 | 1.42857 | 0.98431 |
| 6 | 6 | 0.68515 | 1.44000 | 0.98661 |
| 6 | 7 | 0.70300 | 1.40000 | 0.98419 |
| 6 | 8 | 0.71858 | 1.37143 | 0.98549 |
| 7 | 7 | 0.72647 | 1.36111 | 0.98881 |
| 7 | 8 | 0.74112 | 1.33333 | 0.98817 |
| 8 | 8 | 0.75539 | 1.30612 | 0.98663 |

Table A3. The simulated variance of the within-stratum imbalances in treatment assignments under equal prevalence for 3-factor examples of Pocock and Simon covariate-adaptive allocation with $p=0.9$, 10,000 simulation runs and $500 \times M_s$ subjects per run; the theoretical maximum eigenvalue of the asymptotic correlation matrix of the within-stratum imbalances, and the maximum eigenvalue of the asymptotic covariance matrix of the within-stratum imbalances.

| $n_1$ | $n_2$ | $n_3$ | Simulated variance of the within-stratum imbalances, $\sigma_z^2$ | Theoretical maximum eigenvalue of the asymptotic correlation matrix, $\lambda_{max}$ | Maximum eigenvalue of the asymptotic covariance matrix, $\sigma_z^2 \lambda_{max}$ |
|---|---|---|---|---|---|
| 2 | 2 | 2 | 0.48872 | 2.00000 | 0.97744 |
| 2 | 2 | 3 | 0.57026 | 1.71429 | 0.97759 |
| 2 | 2 | 4 | 0.61348 | 1.60000 | 0.98158 |
| 2 | 2 | 5 | 0.64367 | 1.53846 | 0.99026 |
| 2 | 2 | 6 | 0.65912 | 1.50000 | 0.98868 |
| 2 | 2 | 7 | 0.67175 | 1.47368 | 0.98995 |
| 2 | 2 | 8 | 0.68452 | 1.45455 | 0.99567 |
| 2 | 2 | 9 | 0.68814 | 1.44000 | 0.99092 |
| 2 | 3 | 3 | 0.65497 | 1.50000 | 0.98245 |
| 2 | 3 | 4 | 0.70181 | 1.41176 | 0.99080 |
| 2 | 3 | 5 | 0.72763 | 1.36364 | 0.99223 |
| 2 | 3 | 6 | 0.74207 | 1.33333 | 0.98943 |
| 2 | 3 | 7 | 0.75578 | 1.31250 | 0.99197 |
| 2 | 3 | 8 | 0.76527 | 1.29730 | 0.99278 |
| 2 | 3 | 9 | 0.77156 | 1.28571 | 0.99201 |
| 2 | 4 | 4 | 0.73829 | 1.33333 | 0.98439 |
| 2 | 4 | 5 | 0.76680 | 1.29032 | 0.98941 |
| 2 | 4 | 6 | 0.78504 | 1.26316 | 0.99162 |
| 2 | 4 | 7 | 0.79647 | 1.24444 | 0.99116 |
| 2 | 4 | 8 | 0.80636 | 1.23077 | 0.99244 |
| 2 | 4 | 9 | 0.81105 | 1.22034 | 0.98976 |
| 2 | 5 | 5 | 0.79486 | 1.25000 | 0.99358 |
| 2 | 5 | 6 | 0.81484 | 1.22449 | 0.99776 |
| 2 | 5 | 7 | 0.82580 | 1.20690 | 0.99666 |
| 3 | 3 | 3 | 0.73252 | 1.35000 | 0.98890 |
| 3 | 3 | 4 | 0.77458 | 1.28571 | 0.99588 |
| 3 | 3 | 5 | 0.79431 | 1.25000 | 0.99288 |
| 3 | 3 | 6 | 0.80813 | 1.22727 | 0.99179 |
| 3 | 4 | 4 | 0.80318 | 1.23077 | 0.98854 |
| 3 | 4 | 5 | 0.82471 | 1.20000 | 0.98965 |
| 3 | 4 | 6 | 0.84121 | 1.18033 | 0.99290 |
| 3 | 5 | 5 | 0.84618 | 1.17188 | 0.99161 |
| 3 | 5 | 6 | 0.86126 | 1.15385 | 0.99376 |
| 3 | 6 | 6 | 0.87668 | 1.13684 | 0.99664 |
| 4 | 4 | 4 | 0.83842 | 1.18519 | 0.99368 |
| 5 | 5 | 5 | 0.89275 | 1.11607 | 0.99637 |

Table A4. The simulated variance of the within-stratum imbalances in treatment assignments under equal prevalence for 4- to 7-factor examples of Pocock and Simon covariate-adaptive allocation with $p=0.9$, 10,000 simulation runs and $500 \times M_s$ subjects per run; the theoretical maximum eigenvalue of the asymptotic correlation matrix of the within-stratum imbalances, and the maximum eigenvalue of the asymptotic covariance matrix of the within-stratum imbalances.

| Number of factors $M$ | $n_1$ | $n_2$ | $n_3$ | $n_4$ | $n_5$ | $n_6$ | $n_7$ | Variance of the within-stratum imbalances, $\sigma_z^2$ | Maximum eigenvalue of the asymptotic correlation matrix, $\lambda_{max}$ | Maximum eigenvalue of the asymptotic covariance matrix, $\sigma_z^2 \lambda_{max}$ |
|---|---|---|---|---|---|---|---|---|---|---|
| 4 | 2 | 2 | 2 | 2 | . | . | . | 0.67755 | 1.45455 | 0.98553 |
| 4 | 2 | 2 | 2 | 3 | . | . | . | 0.74462 | 1.33333 | 0.99283 |
| 4 | 2 | 2 | 2 | 4 | . | . | . | 0.77287 | 1.28000 | 0.98927 |
| 4 | 2 | 2 | 2 | 5 | . | . | . | 0.79506 | 1.25000 | 0.99383 |
| 4 | 2 | 2 | 2 | 6 | . | . | . | 0.81248 | 1.23077 | 0.99997 |
| 4 | 2 | 2 | 3 | 3 | . | . | . | 0.79892 | 1.24138 | 0.99177 |
| 4 | 2 | 2 | 3 | 4 | . | . | . | 0.82720 | 1.20000 | 0.99264 |
| 4 | 2 | 2 | 3 | 5 | . | . | . | 0.84606 | 1.17647 | 0.99536 |
| 4 | 2 | 2 | 3 | 6 | . | . | . | 0.85597 | 1.16129 | 0.99403 |
| 4 | 2 | 2 | 4 | 4 | . | . | . | 0.85499 | 1.16364 | 0.99490 |
| 4 | 2 | 2 | 4 | 4 | . | . | . | 0.86904 | 1.14286 | 0.99319 |
| 4 | 2 | 2 | 4 | 4 | . | . | . | 0.88048 | 1.12941 | 0.99442 |
| 5 | 2 | 2 | 2 | 2 | 2 | . | . | 0.81174 | 1.23077 | 0.99907 |
| 5 | 2 | 2 | 2 | 2 | 3 | . | . | 0.84992 | 1.17073 | 0.99503 |
| 6 | 2 | 2 | 2 | 2 | 2 | 2 | . | 0.88907 | 1.12281 | 0.99826 |
| 7 | 2 | 2 | 2 | 2 | 2 | 2 | 2 | 0.93641 | 1.06667 | 0.99884 |

Table A5. The maximum eigenvalue of the asymptotic covariance matrix of the within-stratum imbalances in treatment assignments in examples of the Pocock and Simon covariate-adaptive allocation with two independent factors and $p=0.9$: one factor with two levels of equal prevalence and one factor with three levels of unequal prevalence; 10,000 simulation runs and a sample size of 500,000 per run.

| Prevalence of the three levels of the second factor | Maximum eigenvalue of the estimated asymptotic covariance matrix derived through simulations |
|---|---|
| 1/4, 1/2, 1/4 | 0.98562262 |
| 1/5, 2/5, 2/5 | 0.96420697 |
| 1/5, 1/5, 3/5 | 0.98849695 |
| 1/6, 2/6, 3/6 | 0.99377636 |
| 1/7, 3/7, 3/7 | 0.99504833 |
| 1/7, 2/7, 4/7 | 1.00873787 |

Appendix 2.

Proof of Theorem 2.

To describe the eigenvectors of the matrix $C$ on the set of strata $z=(i_1, \ldots, i_M)$, $i_j=1, \ldots, n_j$, invariant with respect to any permutation of factor levels $i_j=1, \ldots, n_j$ for any $j$, let us use the tensor notations. Instead of labeling the basis vectors in the space where the matrix $C$ acts by $e_{i_1,\ldots,i_M}$, $i_1=1, \ldots, n_1$; ….; $i_M=1, \ldots, n_M$, we will use the tensor product notation $e_{1i_1} \otimes \ldots \otimes e_{Mi_M}$, where $e_{jk}$, $k=1, \ldots, n_j$, is the standard basis in $\mathbf{R}^{n_j}$.

For decomposable vectors $v = v_1 \otimes \ldots \otimes v_M$, where $v_j = \sum_{k=1}^{n_j} v_{jk} e_{jk}$, we have

$$v = \sum_{i_1=1}^{n_1} \ldots \sum_{i_M=1}^{n_M} (v_{1i_1} \ldots v_{Mi_M}) e_{1i_1} \otimes \ldots \otimes e_{Mi_M}$$

Denote $F_j = e_{j1} + \cdots + e_{jn_j}$.

Lemma 2. If the non-negative definite matrix $C$ on the set of strata $z=(i_1, \ldots, i_M)$, $i_j=1, \ldots, n_j$, is invariant with respect to any permutation of factor levels $i_j=1, \ldots, n_j$ for any $j$, that is its terms are

$c((i_1,\ldots,i_M), (j_1,\ldots,j_M)) = \tilde{c}_I, I \subset \{1, \ldots, M\}$, where $I$ denotes the set of the factors common for the two strata $(i_1,\ldots,i_M)$ and $(j_1,\ldots,j_M)$, then

1) The eigenspace that corresponds to the eigenvalue $\lambda_J$ in Theorem 2 is spanned by decomposable vectors $v_1 \otimes \ldots \otimes v_M$ where,
if $j \in J^C$, then $v_j$ is a linear combination of $(e_{jk} - e_{j,k+1})$, $k=1, \ldots, n_j-1$,
and if $j \in J$, then $v_j = F_j$.

2) Eigenvectors do not depend on specific values of $\tilde{c}_I$.

Example. If $M=2$ (2-factor case), the only nonzero eigenvalue is $\lambda_\emptyset = \frac{n_1 n_2}{(n_1-1)(n_2-1)}$.

The corresponding eigenspace is spanned by the vectors

$(e_{1k_1} - e_{1(k_1+1)}) \otimes (e_{2k_2} - e_{2(k_2+1)})$, $k_1=1, \ldots, n_1-1$; $k_2=1, \ldots, n_2-1$.

The rest of the eigenvalues $\lambda_{\{1\}} = \lambda_{\{2\}} = \lambda_{\{1,2\}} = 0$.

The eigenspace corresponding to 0 eigenvalues is spanned by vectors $e_{1k_1} \otimes F_2$, $k_1=1, \ldots, n_1$, and $F_1 \otimes e_{2k_2}$, $k_2=1, \ldots, n_2$

In particular, when both factors have 2 levels ($n_1=n_2=2$), the maximum eigenvalue is 4 and the corresponding eigenvector is (1, -1, -1, 1).

Proof of Lemma 2.

For each $j$, vectors $F_j$ and $(e_{jk} - e_{j,k+1})$, $k=1, \ldots, n_j-1$ form a basis in $\mathbf{R}^{n_j}$ (the last $n_j-1$ vectors span a subspace of vectors whose coordinates add up to 0). Thus, we need to show that vectors of the form $v = v_1 \otimes \ldots \otimes v_M$ described above form an eigenbasis for the correlation matrix $C$.

Let us look at how $C$ acts on the standard basis vectors $e_{1j_1} \otimes \ldots \otimes e_{Mj_M}$ in $V = \mathbf{R}^{n_1} \otimes \ldots \otimes \mathbf{R}^{n_M}$.

In the expansion of $[C](e_{1j_1} \otimes ... \otimes e_{Mj_M})$, a basis vector $e_{1i_1} \otimes ... \otimes e_{Mi_M}$ enters with a coefficient $\tilde{c}_I$ if for $k \in I$, $i_k = j_k$, and, for $k \in I^C$, $i_k \neq j_k$. This can be summed up as

$$[C](e_{1j_1} \otimes ... \otimes e_{Mj_M}) = \sum_{I \subset \{1,...,M\}} \tilde{c}_I u^I, \tag{A1}$$

where $u^I = u_1^I \otimes ... \otimes u_M^I$ has

$$u_k^I = e_{kj_k} \text{ for } k \in I$$

and $\tag{A2}$

$$u_k^I = \sum_{\substack{i=1 \\ i \neq j_k}}^{n_k} e_{ki} = F_k - e_{kj_k} \text{ for } k \in I^C.$$

Next, fix $J$ and let $v^J = v_1 \otimes ... \otimes v_M$ be a vector in $V$ such that

for $k \in J$, $v_k = F_k$ and

for $k \in J^C$, $v_k$ is a vector in $R^{n_k}$ whose coordinates add up to 0, that is $v_k = \sum_{j=1}^{n_k} \alpha_{kj} e_{kj}$, $\sum_{j=1}^{n_k} \alpha_{kj} = 0$.

Using (A1) and (A2) and the multilinearity of the tensor product, we see that for each $I$, the components of $u^I$ in (A1) add up to a vector in $V$ whose $k$-th tensor factor is:

$F_k = v_k$ for $k \in I \cap J$;

$\sum_{j=1}^{n_k} \alpha_{kj}(F_k - e_{kj}) = \left(\sum_{j=1}^{n_k} \alpha_{kj}\right) F_k - \sum_{j=1}^{n_k} \alpha_{kj} e_{kj} = -v_k$ for $k \in I^C \cap J^C$;

$\sum_{j=1}^{n_k}(F_k - e_{kj}) = (n_k - 1) F_k = (n_k - 1) v_k \qquad$ for $k \in I^C \cap J$;

and $\sum_{j=1}^{n_k} \alpha_{kj} e_{kj} = v_k$ $k \in I \cap J^C$.

In all four cases we obtained a vector proportional to $v_k$. This means

$$Cv^J = \left(\sum_{I \subset \{1,...M\}} (-1)^{\#J^C \cap I^C} \times \prod_{k \in J \cap I^C} (n_k - 1)\right) v^J$$

which proves that

$$\lambda_J = \sum_{I \subset \{1,...M\}} (-1)^{\#J^C \cap I^C} \times \tilde{c}_I \times \prod_{j \in J \cap I^C} (n_j - 1)$$

described by (37) is the eigenvalue of the matrix $C$.

It remains to compute the multiplicity i.e. the dimension of the space spanned by the vectors of the form $v^J$ which is equal to the product of dimensions of "zero-sum" subspaces in $R^{n_k}$ for $k \in J^C$, that is $m_J = \prod_{j \in J^C} (n_j - 1)$. This completes the proof of Theorem 2.

Proof of Lemma 1.

Let us denote by $Q$ the denominator of (35): $Q = \prod_{i=1}^M n_i - \sum_{i=1}^M n_i + M - 1$.

Let us denote by $q_i = n_i - 1, i = 1, ..., M$. Then $Q$ can be written as

$$Q = \prod_{i=1}^M (q_i + 1) - \sum_{i=1}^M q_i - 1 \tag{A3}$$

And (35) can be written as

$$c_I = \frac{\#I^C - 1 - \sum_{i \in I} q_i}{Q}, \text{ when } I \neq \{1, \ldots, M\};$$  (A4)

$c_I = 1$, when $I = \{1, \ldots, M\}$.

Consider a case when $\#J^C = 0$ that is $J = \{1, \ldots, M\}$; let us prove that $\lambda_J = 0$.

Then the expression (37) for eigenvalue $\lambda_J$ and $\tilde{c}_I = c_I$ can be written as

$$\lambda_J = \sum_{I \subset \{1,\ldots M\}} c_I \times \prod_{j \in I^C} q_j = \frac{1}{Q}\left(Q + \sum_{I \subsetneq \{1,\ldots M\}}(\#I^C - 1 - \sum_{i \in I} q_i) \times \prod_{j \in I^C} q_j\right)$$  (A5)

Denote the expression in parentheses in (A5) by

$X = Q + \sum_{I \subsetneq \{1,\ldots M\}}(\#I^C - 1 - \sum_{i \in I} q_i) \times \prod_{j \in I^C} q_j$; $X$ is a polynomial in $q_1, \ldots, q_M$. Let us denote by

$Y = \sum_{I \subsetneq \{1,\ldots M\}}(\#I^C - 1 - \sum_{i \in I} q_i) \times \prod_{j \in I^C} q_j$; $X = Q + Y$.

Let us compute the coefficients of the polynomial $X$. From (A3), $Q$ has zero constant term and zero linear terms. Every term of $Y$ contains at least one of the factors $q_1, \ldots, q_M$, so $Y$ also has zero constant term. For $Y$ to have a linear term $q_i$, requires $I^C = q_i$, but then the coefficient of $q_i$ is $\#I^C - 1 = 0$. Thus, $Y$ also has only terms of degree 2 or higher.

Now consider any monomial $q_{j_1} \ldots q_{j_k}, k > 1$. It enters $Q$ with the coefficient 1. In $Y$, this monomial appears when $I^C = \{j_1, \ldots, j_k\}$ with the coefficient $\#I^C - 1 = k - 1$, or when $I^C = \{j_1, \ldots, j_k\} \setminus \{j_s\}$, $s = 1, \ldots k$, and $q_{j_s}$ enters the sum $\sum_{i \in I} q_i$ in (A4) in which case the corresponding coefficient is $(-1)$. Therefore, in total, the monomial $q_{j_1} \ldots q_{j_k}, k > 1$ appears in $Y$ with the coefficient $(k - 1) - k = -1$ and thus cancels with the corresponding term in $Q$. This proves that $X = 0$ and thus $\lambda_J = 0$.

Now consider a case when $\#J^C = 1$; let us prove that $\lambda_J = 0$. Due to symmetry with respect to factors, it is sufficient to consider the case of $J = \{2, \ldots, M\}$. In this case from (37) with $\tilde{c}_I = c_I$,

$$\lambda_J = \frac{1}{Q}\left(Q + \sum_{I \subsetneq \{1,\ldots M\}}(\#I^C - 1 - \sum_{i \in I} q_i) \times \prod_{j \in I^C} q_j - \sum_{I \subset \{2,\ldots M\}}(\#I^C - 1 - \sum_{i \in I} q_i) \times \prod_{j \in I^C \setminus \{1\}} q_j\right)$$  (A6)

Let us denote $Y_1 = \sum_{I \subsetneq \{1,\ldots M\}}(\#I^C - 1 - \sum_{i \in I} q_i) \times \prod_{j \in I^C} q_j$ and

$Y_2 = \sum_{I \subset \{2,\ldots M\}}(\#I^C - 1 - \sum_{i \in I} q_i) \times \prod_{j \in I^C \setminus \{1\}} q_j$.

Let us prove that the expression in parenthesis of (A6) $X_2 = Q + Y_1 - Y_2$ is identically 0 as a polynomial in $q_1, \ldots, q_M$. In $Y_2$, the constant term corresponds to $I = \{2, \ldots, M\}$ and appears with the coefficient $\#I^C - 1 = 0$. A monomial $q_{j_1} \ldots q_{j_k}, k \geq 1$; $j_1, \ldots, j_k \in \{2, \ldots, M\}$ appear in $Y_2$ either when $I^C = \{1, j_1, \ldots, j_k\}$ or $I^C = \{1, j_1, \ldots, j_k\} \setminus \{j_s\}, s = 1, \ldots k$, with the coefficient in total $(\#I^C - 1) - k = 0$. Thus $Y_2 = 0$.

The fact that $Q + Y_1 = 0$ can be proven using the argument analogous to the argument used above for the case $J = \{1, \ldots, M\}$. Thus, $X_2 = 0$ and $\lambda_J = 0$ for $J = \{1, \ldots, M\}$.

Now to prove that $\lambda_J = \lambda_{max}$ described by (38) for all $J$ such that $\#J < M - 1$, we will use the following identities that are a restatement of well-known properties of binomial coefficients:

$\sum_{I \subset S}(-1)^{\#I} = 0$ (#$S \geq 1$)

(A7)

$\sum_{I \subset S}(-1)^{\#I} \#I = 0$ (#$S \geq 2$).

Let us denote $r_i = 1 - n_i, i = 1, ..., M$ and let us re-write $Q$, $c_I$, and $\lambda_J$ through $r_i$:

$Q = \prod_{i=1}^{M}(1 - r_i) + \sum_{i=1}^{M} r_i - 1$

$c_I = \frac{\#I^C - 1 + \sum_{i \in I} r_i}{Q}$, when $I \neq \{1, ..., M\}$;

$\lambda_J = 1 + \sum_{I \subsetneq \{1,...M\}}(-1)^{\#I^C} \frac{\#I^C - 1 + \sum_{i \in I} r_i}{Q} \times \prod_{j \in I^C} r_j$ (A8)

The sum in (A8) is a ratio with the denominator $Q$ and the numerator denoted $X_3$ which is the polynomial in $r_1, ..., r_M$.

First, we will show that $X_3$ is at most linear in $r_1, ..., r_M$. Indeed, $X_3$ can only contain nonlinear monomials of the form

$r_{j_1} ... r_{j_k}, k \geq 1; \{j_1, ..., j_k\} \subset J, k \geq 2$

or (A9)

$r_{j_1} ... r_{j_k} r_i, \{j_1, ..., j_k\} \subset J, i \in J^C, k \geq 1$

The latter can only appear if $J \cap I^C = \{j_1, ..., j_k\}$ and $i \in I$, that is $I = (J \setminus \{j_1, ..., j_k\}) \cup \{i\} \cup I'$ where $I'$ is an arbitrary subset of $J^C$. Thus the coefficient of $r_{j_1} ... r_{j_k} r_i$ in $X_3$, up to a common sign, is

$\sum_{I' \subset J^C}(-1)^{\#(I')} = 0$ from (A7).

The first monomial in (A9) appears in $X_3$ when $J \cap I^C = \{j_1, ..., j_k\}$ or $J \cap I^C = \{j_1, ..., j_k\} \setminus \{j_s\}$ for some $s \in \{1, ..., k\}$. The corresponding coefficient is then (again, up to a common sign)

$\sum_{(I')^C \subset J^C}(-1)^{\#(I')^C}(\#(I')^C - 1 - k) = 0$ from (A7).

We conclude that $X_3$ does not contain monomials in $r_1, ..., r_M$ of degree 2 or higher.

Next, a linear monomial $r_i$ enters $X_3$ only $J \cap I^C = \emptyset$ and $i \in I$, that is $I = (J \cup \{i\}) \cup I'$, where $I' \subsetneq (J \cup \{i\})^C$. The corresponding coefficient is

$\sum_{\emptyset \subsetneq I^C \subset (J \cup \{i\})^C}(-1)^{\#I^C} = -1$ (the sum contains one less term than the identity in (A7).

Similarly, the constant term in $X_3$ corresponds to a subset I such that $J \subset I \subsetneq \{1, ..., M\}$ and has a coefficient $\sum_{\emptyset \subsetneq I^C \subset J^C}(-1)^{\#I^C}(\#I^C - 1) = 1$.

To summarize,

$X_3 = -\sum_{i=1}^{M} r_i + 1$

and

$$\lambda_J = 1 + \frac{1 - r_1 - \ldots - r_M}{Q} = \frac{\prod_{i=1}^{M}(1-r_i)}{Q} = \frac{\prod_{i=1}^{M} n_i}{Q}.$$

This proves (38) and completes the proof of Lemma 1.